\documentclass[11pt,reqno]{amsart}
\usepackage[margin=1in]{geometry}
\usepackage{amsmath,amssymb,amsthm}
\usepackage{tikz}
\usetikzlibrary{arrows.meta}
\usepackage[colorlinks=true,linkcolor=black,citecolor=black,urlcolor=blue]{hyperref}

\theoremstyle{plain}
\newtheorem{theorem}{Theorem}[section]
\newtheorem{proposition}[theorem]{Proposition}
\newtheorem{lemma}[theorem]{Lemma}
\newtheorem{corollary}[theorem]{Corollary}
\theoremstyle{remark}
\newtheorem{remark}[theorem]{Remark}

\numberwithin{equation}{section}

\title{Fourier Ratio Bounds and Sensitivity to Sensor Placement}

\author[N. Alberti]{Nico Alberti}
\address{Department of Mathematics, Virginia Tech, Blacksburg, VA 24061}
\email{nicoalberti@vt.edu}

\author[C. Housholder]{Christopher Housholder}
\address{Department of Mathematics, Missouri State University, Springfield, MO 65897}
\email{christopherlancehousholder@gmail.com}

\author[A. Iosevich]{Alex Iosevich}
\address{Department of Mathematics, University of Rochester, Rochester, NY 14627}
\email{alex.iosevich@rochester.edu}
\thanks{The third author was supported in part by NSF grant DMS-2506858 and a Simons Fellowship in Mathematics in 2026.}

\author[S. J. Miller]{Steven J. Miller}
\address{Department of Mathematics, Williams College, Williamstown, MA 01267}
\email{sjm1@williams.edu}

\author[E. A. Palsson]{Eyvindur A. Palsson}
\address{Department of Mathematics, Virginia Tech, Blacksburg, VA 24061}
\email{palsson@vt.edu}

\begin{document}

\begin{abstract}
We study how the regularity of a continuous function controls the Fourier ratio of its sampled values, and how errors in sensor placement affect recovery from incomplete measurements. For a periodic function on the square sampled on an $N$ by $N$ grid, we show that H\"older regularity $C^{0,\gamma}$ with $0<\gamma<1$ leads to a bound of order $N^{1-\gamma}$, that Lipschitz regularity leads to the critical bound $\sqrt{1+\log N}$, and that fractional Sobolev regularity $W^{s,p}$ with $0<s<1$ and $sp>2$ leads to a bound of order $N^{1-s}$, where both polynomial exponents are sharp. A common translation of the grid removes the hypothesis on the sampled mass needed on a fixed grid and yields an estimate in terms of the $L^2$ norm of the gradient, and with two derivatives the bound is independent of the grid size. We extend the energy argument behind these results to compact Riemannian manifolds and prove that the critical square-root logarithm is necessary on every compact surface. Finally, we treat separate sensor displacements as measurement error and combine the analytic bounds with random sampling to obtain stable recovery of the sampled vector on the square and of bandlimited functions on compact manifolds, where the number of samples depends on the suprema of the chosen eigenfunctions.
\end{abstract}

\maketitle

\tableofcontents

\section{Introduction}

The purpose of this paper is to study how the regularity of a continuous function controls the concentration of the Fourier coefficients of its samples, and what happens to the recovery of the function from incomplete measurements when the sensors that take those samples are moved. The recovery of a signal from incomplete measurements has been studied extensively in harmonic analysis, sampling theory, compressed sensing, and inverse problems (see, for example, \cite{BIN,CRTRobust,CRTStable,FoucartRauhut,IPY,Rauhut,RV}), and the underlying principle is that if the Fourier coefficients of a signal are sufficiently concentrated, then stable recovery may be possible from a small portion of its spatial samples.

A convenient way to measure this concentration is the Fourier ratio \cite{AldalehEtAl,BIN}. Let $\mathbb Z_N=\mathbb Z/N\mathbb Z$, where $N\geq2$ is an integer, and given $g:\mathbb Z_N^d\to\mathbb C$ define the Fourier transform of $g$ by
$$ \widehat g(m) = \frac{1}{N^{\frac{d}{2}}} \sum_{x\in\mathbb Z_N^d} e^{-\frac{2\pi i x\cdot m}{N}}g(x). $$
With this normalization the Fourier transform is unitary, which is to say that $\|\widehat g\|_2=\|g\|_2$, where $\|\cdot\|_1$ and $\|\cdot\|_2$ denote the $\ell^1$ and $\ell^2$ norms on $\mathbb Z_N^d$. If $g$ is not identically zero, its Fourier ratio is defined by
$$ FR(g)=\frac{\|\widehat g\|_1}{\|\widehat g\|_2}. $$
If precisely $k$ Fourier coefficients of $g$ are nonzero, then $FR(g)\leq\sqrt{k}$ by Cauchy--Schwarz, but the ratio remains meaningful when many coefficients are nonzero and most of them are small. On an $N$ by $N$ grid one always has $1\leq FR(g)\leq N$, so a Fourier ratio substantially smaller than $N$ indicates spectral concentration, and existing recovery theorems turn such a bound into a guarantee for reconstruction from incomplete random measurements \cite{BIN,CRTRobust,FoucartRauhut,IPY}.

Suppose now that the discrete data comes from sensors that sample a continuous function. It is natural to ask what the regularity of this function tells us about its sampled Fourier coefficients, and what happens if the sensors are moved. We shall see that H\"older and Sobolev regularity control the Fourier ratio on finite grids, that averaging over a common translation of the grid strengthens the estimates, and that separate placement errors can be included in the recovery argument as measurement noise.

Our starting point is the work of Iosevich, Palsson, and Yavicoli \cite{IPY}, who studied discretization and sampling through the Fourier ratio. For a $C^2$ periodic function sampled on an $N$ by $N$ grid they obtained a bound of order $\log N$ and used it in a recovery theorem for random samples on $\mathbb Z_N^2$, and they also considered bandlimited functions on $S^2$ by means of quadrature and spherical harmonic expansions. The Fourier ratio has also been connected to recovery, localization, learning, PDE propagation, spectral synthesis, and uncertainty principles in \cite{AldalehEtAl,BIN,IIPY,IMWUncertainty,IMW,IP}. We keep the point of view of \cite{IPY}, but in place of separate pointwise estimates for the Fourier coefficients we use an energy argument.

Our first result shows that differentiability is not required. Let $f$ be a periodic function on $[0,1]^2$ and let $g(x)=f(\frac{x}{N})$, $x\in\mathbb Z_N^2$, denote its samples on the $N$ by $N$ grid. If $f$ belongs to the H\"older class $C^{0,\gamma}([0,1]^2)$ with $0<\gamma<1$, then for sufficiently fine grids
$$ FR(g)\leq C_\gamma \frac{\|f\|_{C^{0,\gamma}}}{\|f\|_{L^2}}N^{1-\gamma}, $$
where the H\"older norm, like the other norms that appear in this introduction, is defined in Section~\ref{sec:square}. In particular, every positive H\"older exponent improves on the universal bound $FR(g)\leq N$. At the endpoint $\gamma=1$ we obtain
$$ FR(g)\leq C \frac{\|f\|_{C^{0,1}}}{\|f\|_{L^2}}\sqrt{1+\log N}, $$
which means that Lipschitz regularity already reaches the critical estimate for one derivative. For $C^2$ functions an argument with second differences yields a bound independent of $N$ on sufficiently fine fixed grids, or on a suitable translated grid for every $N\geq2$. The energy calculation thus lowers the regularity needed for a nontrivial estimate, and it also improves the bound at the regularity considered in \cite{IPY}.

The H\"older result has a fractional Sobolev counterpart. Suppose that $0<s<1$, $1<p<\infty$, and $sp>2$. A local form of the Morrey inequality, used inside the same dyadic argument, shows that for sufficiently fine grids
$$ FR(g)\leq C_{s,p} \frac{\|f\|_{W^{s,p}}}{\|f\|_{L^2}} N^{1-s}. $$
The condition $sp>2$ is there to provide the continuous representative needed for point sampling, and it does not affect the exponent. A sampled discrete chirp shows that neither $1-\gamma$ nor $1-s$ can be lowered in bounds that are uniform over the corresponding regularity classes.

Let us describe why energy enters the picture. A finite difference of the sampled signal multiplies its Fourier coefficients by a discrete frequency, so Parseval's identity controls a weighted $\ell^2$ sum of those coefficients, and Cauchy--Schwarz then reduces the $\ell^1$ estimate to a sum of inverse frequency weights. In dimension two, one derivative leaves a square-root logarithm. Below the Lipschitz endpoint we apply differences at the scale of each frequency annulus, which avoids an unnecessary logarithmic loss, and above the critical number of derivatives the inverse weights are summable.

The estimates on a fixed grid require a hypothesis which prevents the sampled $L^2$ mass from becoming too small. For a $C^1$ function a common translation of the grid removes this issue. Averaging the sampled mass and the difference energy together produces, for every $N\geq2$, a translation $t$ such that the translated samples $g_t(x)=f(\frac{x+t}{N})$ satisfy
$$ FR(g_t)\leq 1+\frac{1}{\sqrt{2}}\sqrt{1+\log N} \frac{\|\nabla f\|_{L^2([0,1]^2)}}{\|f\|_{L^2([0,1]^2)}}. $$
The translation may depend on $f$, so this is an existence theorem for a favorable grid and not a procedure for finding that grid from incomplete observations.

\begin{figure}[t]
\centering
\begin{tikzpicture}[scale=0.82,>=Stealth]
\begin{scope}
\draw[step=1,gray!55,thin] (0,0) grid (3,3);
\foreach \x in {0,1,2}
  \foreach \y in {0,1,2}{
    \draw[->,gray!75,thin] (\x,\y) -- ({\x+0.28},{\y+0.32});
    \draw[fill=white] (\x,\y) circle (1.7pt);
    \fill[blue!70!black] ({\x+0.28},{\y+0.32}) circle (1.7pt);
  }
\node at (1.5,-0.55) {(a) One common translation};
\end{scope}
\begin{scope}[xshift=7.2cm]
\draw[step=1,gray!55,thin] (0,0) grid (3,3);
\foreach \x/\y/\a/\b in {
0/0/0.18/0.22,
1/0/0.45/0.16,
2/0/0.25/0.52,
0/1/0.40/0.28,
1/1/0.20/0.48,
2/1/0.52/0.24,
0/2/0.22/0.42,
1/2/0.48/0.30,
2/2/0.32/0.18}{
    \draw[->,gray!75,thin] (\x,\y) -- ({\x+\a},{\y+\b});
    \draw[fill=white] (\x,\y) circle (1.7pt);
    \fill[red!70!black] ({\x+\a},{\y+\b}) circle (1.7pt);
  }
\node at (1.5,-0.55) {(b) Separate perturbations};
\end{scope}
\end{tikzpicture}
\caption{The open circles are the reference grid points. In (a) every point is moved by the same vector, so the translated points still form a grid, while in (b) each sensor is displaced separately inside its cell and the displacements are treated as measurement errors.}
\label{fig:grid-perturbations}
\end{figure}
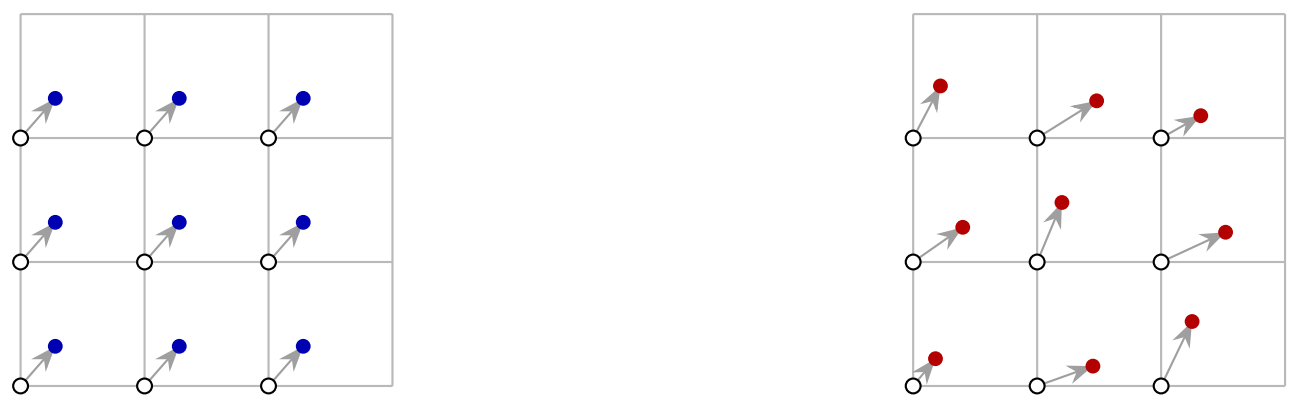

A common translation and a separate error in each sensor position play different roles, as illustrated in Figure~\ref{fig:grid-perturbations}. A common translation preserves the grid, while separate displacements are treated as measurement noise. We prove a recovery theorem for random grid indices and arbitrary displacements inside the translated cells, with an error controlled by the Fourier ratio and by $\frac{1}{N}\|f\|_{C^1}$. If the positions are chosen uniformly inside these cells, the result is a collection of independent uniform sensors on the torus. The reconstruction theorem concerns the sampled vector, and a piecewise constant extension leads to an $L^2$ estimate on the torus with the same placement scale.

The discrete calculation has a spectral analogue. On a connected compact Riemannian manifold of dimension $d$, weighted Cauchy--Schwarz and Weyl's law lead to the critical Sobolev index $s=\frac{d}{2}$. We use $\Lambda$ for an eigenvalue cutoff and $L$ for a spherical harmonic degree. On a compact surface, coefficients proportional to $\frac{1}{\lambda_j}$, where the $\lambda_j$ are the positive eigenvalues of the Laplace--Beltrami operator, show that the factor $\sqrt{1+\log\Lambda}$ is necessary in the estimate involving the $L^2$ norm of the gradient. This establishes sharpness for the critical Sobolev energy, but not for estimates involving the stronger $C^1$ norm.

For recovery we use the bounded orthonormal system theorem of Brugiapaglia, Dirksen, Jung, and Rauhut \cite{BDJR}, together with the stable $\ell^1$ theory in \cite{CRTStable,FoucartRauhut}, and we refer to \cite{BIN} for the formulation in terms of the Fourier ratio. On a manifold the $L^\infty$ bound for the eigenfunctions enters the sample count. Moreover, when an eigenvalue has multiplicity, the $\ell^1$ norm of the coefficients depends on the basis chosen in that eigenspace, and for this reason all the statements on manifolds are made relative to one fixed orthonormal eigenbasis.

This paper is organized as follows. In Section~\ref{sec:tools} we prove the discrete energy estimate and state the bounded orthonormal system theorem. Section~\ref{sec:square} contains the estimates on the square, first for H\"older and fractional Sobolev regularity and then for translated grids and two derivatives. In Section~\ref{sec:manifolds} we treat compact manifolds and the sharp critical estimate on surfaces, and in Section~\ref{sec:recovery} we combine these bounds with sensor placement estimates and stable reconstruction on the square and on compact manifolds.

Throughout the paper the grid size $N\geq2$, the dimension $d\geq1$, and all sample counts are integers, and the error parameter $\varepsilon$ is positive. We write $A\asymp B$ if $cB\leq A\leq CB$ with positive constants $c$ and $C$ independent of the grid size or the spectral cutoff, and subscripts on constants indicate what they are allowed to depend on. We also write $\|\cdot\|_{\ell^2}$ in place of $\|\cdot\|_2$ when it is helpful to distinguish the discrete norm from an $L^2$ norm.

\section{Tools}\label{sec:tools}

We begin with the two tools that are used throughout the paper. The discrete energy estimate controls the Fourier ratio, and the bounded orthonormal system theorem turns this control into recovery from random measurements.

\subsection{Discrete Fourier energy}

The argument that we use on the square in Section~\ref{sec:square} is not special to two dimensions, so we first prove the underlying statement on the finite group $\mathbb Z_N^d$, which also makes it clear why dimension two is critical for one derivative. Let $e_1,\ldots,e_d$ be the standard coordinate vectors, and for $m\in\mathbb Z_N^d$ define $\rho(m)\geq0$ by
$$ \rho(m)^2 = \sum_{j=1}^d \left|e^{\frac{2\pi i m_j}{N}}-1\right|^2. $$
The quantity $\rho(m)$ is zero only at $m=0$, and it plays the role of the size of the frequency $m$ on the finite group.

\begin{proposition}\label{prop:discrete-energy}
Let $g:\mathbb Z_N^d\to\mathbb C$ be nonzero, and let $s>0$. Then
$$ FR(g) \leq \frac{|\widehat g(0)|}{\|g\|_2} + \left(\sum_{m\neq0}\frac{1}{\rho(m)^{2s}}\right)^{\frac{1}{2}} \frac{ \left(\sum_m\rho(m)^{2s}|\widehat g(m)|^2\right)^{\frac{1}{2}} }{\|g\|_2}. $$
Moreover,
$$ \left(\sum_{m\neq0}\frac{1}{\rho(m)^{2s}}\right)^{\frac{1}{2}} \leq C_{d,s} \begin{cases} N^s, & s>\frac{d}{2},\\ N^{\frac{d}{2}}\sqrt{1+\log N}, & s=\frac{d}{2},\\ N^{\frac{d}{2}}, & 0<s<\frac{d}{2}, \end{cases} $$
and when $d=2$ and $s=1$ one may take the more explicit estimate
$$ \sum_{m\neq0}\frac{1}{\rho(m)^2} \leq \frac{N^2}{2}(1+\log N). $$
Finally, if $s=1$ and $\Delta_jg(x)=g(x+e_j)-g(x)$, then
$$ \sum_m\rho(m)^2|\widehat g(m)|^2 = \sum_{j=1}^d\|\Delta_jg\|_2^2. $$
\end{proposition}

\begin{proof}
We begin with the energy identity. A change of variables shows that
$$ \widehat{\Delta_jg}(m) = \frac{1}{N^{\frac{d}{2}}} \sum_x e^{-\frac{2\pi i x\cdot m}{N}}g(x+e_j)-\widehat g(m) = \left(e^{\frac{2\pi i m_j}{N}}-1\right)\widehat g(m), $$
so by Parseval's identity
$$ \|\Delta_jg\|_2^2 = \sum_m \left|e^{\frac{2\pi i m_j}{N}}-1\right|^2 |\widehat g(m)|^2, $$
and summing over $j$ we obtain the last identity in the statement.

To prove the main estimate we separate the zero frequency and apply Cauchy--Schwarz to the remaining frequencies, which leads to
$$ \sum_{m\neq0}|\widehat g(m)| = \sum_{m\neq0} \frac{1}{\rho(m)^s}\rho(m)^s|\widehat g(m)| \leq \left(\sum_{m\neq0}\frac{1}{\rho(m)^{2s}}\right)^{\frac{1}{2}} \left(\sum_m\rho(m)^{2s}|\widehat g(m)|^2\right)^{\frac{1}{2}}. $$
Adding $|\widehat g(0)|$ and dividing by $\|\widehat g\|_2=\|g\|_2$ we arrive at the first assertion.

It remains to estimate the first factor. Represent an element $k\in\mathbb Z_N$ by an integer in $\{0,\ldots,N-1\}$ and write $|k|_*=\min\{k,N-k\}$. Since $\sin u\geq\frac{2u}{\pi}$ for $0\leq u\leq\frac{\pi}{2}$, we have
$$ \left|e^{\frac{2\pi i k}{N}}-1\right| = 2\sin\left(\frac{\pi|k|_*}{N}\right) \geq \frac{4|k|_*}{N}, $$
and it follows that
$$ \rho(m)^2 \geq \frac{16}{N^2}\sum_{j=1}^d|m_j|_*^2. $$
We group the nonzero frequencies according to the value of $r=\max_{1\leq j\leq d}|m_j|_*$. The number of frequencies on the shell indexed by $r$ is at most $(2r+1)^d-(2r-1)^d\leq C_dr^{d-1}$, and every frequency on that shell satisfies $\sum_{j=1}^d|m_j|_*^2\geq r^2$. Consequently,
$$ \sum_{m\neq0}\frac{1}{\rho(m)^{2s}} \leq C_{d,s}N^{2s} \sum_{r=1}^{\left\lfloor\frac{N}{2}\right\rfloor}r^{d-1-2s}. $$
The sum on the right is bounded independently of $N$ when $2s>d$, by $C(1+\log N)$ when $2s=d$, and by $C_{d,s}N^{d-2s}$ when $2s<d$, and the three cases follow upon taking square roots. When $d=2$ and $s=1$, the shell indexed by $r$ has at most $8r$ points, and the same calculation shows that
$$ \sum_{m\neq0}\frac{1}{\rho(m)^2} \leq \frac{N^2}{16} \sum_{r=1}^{\left\lfloor\frac{N}{2}\right\rfloor}\frac{8r}{r^2} \leq \frac{N^2}{2}(1+\log N), $$
as claimed.
\end{proof}

The three regimes have a simple interpretation. When a smooth function is sampled, a discrete derivative of order $s$ contributes a factor of approximately $\frac{1}{N^s}$. Once this scaling is taken into account, the Fourier ratio is bounded independently of $N$ above the critical index, grows like $\sqrt{\log N}$ at the critical index, and grows like $N^{\frac{d}{2}-s}$ below it. In Section~\ref{sec:square} we push the argument below one derivative using H\"older increments and the fractional Sobolev--Morrey embedding. The endpoint $\gamma=1$ is the critical case in two dimensions, and it already contains the $C^1$ functions.

\subsection{Bounded orthonormal system recovery}

We now state the standard recovery result that we shall use. Let $a=(a_1,\ldots,a_D)\in\mathbb C^D$. Given an integer $s\geq1$, a vector is called $s$-sparse if at most $s$ of its coordinates are nonzero, and for $1\leq s\leq D$ the best $s$-term approximation error of $a$ in $\ell^1$ is
$$ \sigma_s(a)_1 = \inf\{\|a-b\|_1: b\text{ has at most }s\text{ nonzero coordinates}\}. $$
A minimizer is obtained by retaining the $s$ largest coordinates of $a$ in absolute value, and the quantity $\sigma_s(a)_1$ measures how close $a$ is to being $s$-sparse.

Let $\mu$ be a probability measure on a set $\Omega$, and suppose that $\phi_1,\ldots,\phi_D$ are orthonormal in $L^2(\Omega,\mu)$ and satisfy $\|\phi_j\|_{L^\infty(\Omega)}\leq K$, where $K\geq1$ necessarily by orthonormality. Choose $x_1,\ldots,x_q$ independently according to $\mu$ and form the sampling matrix
$$ A_{ij}=\frac{1}{\sqrt q}\phi_j(x_i), $$
where the factor $\frac{1}{\sqrt q}$ makes the expected squared norm of $Aa$ equal to $\|a\|_2^2$. The matrix $A$ is said to have the restricted isometry property of order $s$ with constant $\delta$ if $0<\delta<1$ and
$$ (1-\delta)\|a\|_2^2 \leq \|Aa\|_2^2 \leq (1+\delta)\|a\|_2^2 $$
for every $s$-sparse vector $a$. In other words, $A$ acts almost like an isometry on every coordinate subspace of dimension at most $s$.

\begin{proposition}[Bounded orthonormal system recovery]\label{prop:bos-recovery}
Let $\phi_1,\ldots,\phi_D$ and $A$ be as above. Fix $1\leq s\leq D$ and $\gamma\geq1$. There is a constant $C_\gamma$ such that, if
$$ q \geq C_\gamma K^2s\log^2(2sK^2)\log(2D), $$
then, with probability at least $1-\frac{1}{D^\gamma}$, the following holds simultaneously for all coefficient vectors $a$ and all noise vectors $e$. If $\eta\geq0$, $y=Aa+e$ with $\|e\|_2\leq\eta$, and $a^*$ minimizes $\|b\|_1$ subject to $\|Ab-y\|_2\leq\eta$, then
$$ \|a^*-a\|_2 \leq 3\frac{\sigma_s(a)_1}{\sqrt{s}}+5\eta. $$
\end{proposition}

\begin{proof}
Let $k=\min\{2s,D\}$ and $\delta=\frac{1}{16}$. By Theorem~2.3 of Brugiapaglia, Dirksen, Jung, and Rauhut \cite{BDJR}, applied to the orthonormal system above, the matrix $A$ has the restricted isometry property of order $k$ with constant $\delta$ provided that
$$ q\geq C K^2k\log^2(2kK^2)\log(2D), $$
with failure probability at most
$$ 2\exp\left(-\frac{cq}{K^2k}\right), $$
where $C$ and $c>0$ are absolute constants because $\delta$ is fixed. Since $k\leq2s$ and $K\geq1$, the hypothesis of the proposition implies both this bound on the sample size and the required failure probability after increasing $C_\gamma$.

We include the recovery argument in order to keep track of the constants. We work on the restricted isometry event and write $h=a^*-a$. Both $a$ and $a^*$ satisfy the measurement constraint, so $\|Ah\|_2\leq2\eta$. If $2s>D$, restricted isometry holds on all of $\mathbb C^D$, hence
$$ \|h\|_2\leq\frac{2\eta}{\sqrt{1-\delta}}\leq3\eta, $$
and this includes the case $s=D$.

Suppose now that $2s\leq D$, and let $S$ index the $s$ largest coordinates of $a$ in absolute value. For a subset $E$ of the coordinate index set we denote by $h_E$ the vector equal to $h$ on $E$ and zero elsewhere, and by $E^c$ the complement of $E$. The minimizing property of $a^*$ and the triangle inequality imply that
$$ \|h_{S^c}\|_1\leq\|h_S\|_1+2\sigma_s(a)_1. $$
Partition $S^c$ into sets $S_1,S_2,\ldots$ of size $s$, except possibly the last, in decreasing order of $|h_j|$. Every coordinate on $S_j$ is at most $\frac{1}{s}\|h_{S_{j-1}}\|_1$ when $j\geq2$, and therefore
$$ \sum_{j\geq2}\|h_{S_j}\|_2 \leq\frac{\|h_{S^c}\|_1}{\sqrt{s}} \leq\|h_S\|_2+2\frac{\sigma_s(a)_1}{\sqrt{s}}. $$
With $T=S\cup S_1$, $H=\|h_T\|_2$, and $u=\frac{\sigma_s(a)_1}{\sqrt{s}}$, this estimate implies that
$$ \sum_{j\geq2}\|h_{S_j}\|_2\leq H+2u \text{ and } \|h\|_2\leq2H+2u. $$

If $v$ and $w$ are disjointly supported vectors, each with at most $s$ nonzero coordinates, then by restricted isometry
$$ |\langle Av,Aw\rangle|\leq\delta\|v\|_2\|w\|_2, $$
because on the union of their supports the operator $A^*A-I$ has norm at most $\delta$. Applying this estimate separately to $h_S$ and $h_{S_1}$ we obtain
$$ |\langle Ah_T,Ah_{S_j}\rangle| \leq\sqrt{2}\delta H\|h_{S_j}\|_2 \text{ for } j\geq2, $$
and since $Ah_T=Ah-\sum_{j\geq2}Ah_{S_j}$, it follows that
$$ (1-\delta)H^2 \leq2\sqrt{1+\delta}\eta H +\sqrt{2}\delta H\sum_{j\geq2}\|h_{S_j}\|_2. $$
If $H=0$, the ordering of the blocks implies that $h=0$. Otherwise we divide by $(1-\delta)H$ and use the tail bound. Writing
$$ \alpha=\frac{2\sqrt{1+\delta}}{1-\delta} \text{ and } \beta=\frac{\sqrt{2}\delta}{1-\delta}, $$
we have $H\leq\alpha\eta+\beta(H+2u)$, and consequently
$$ \|h\|_2 \leq\frac{2(1+\beta)}{1-\beta}u +\frac{2\alpha}{1-\beta}\eta. $$
When $\delta=\frac{1}{16}$ we have $\beta<\frac{1}{10}$ and $\alpha<\frac{9}{4}$, so the two coefficients on the right are less than $3$ and $5$, respectively, and the proof is complete.
\end{proof}

The sampling input in this proof is the bounded orthonormal system case of \cite[Theorem~2.3]{BDJR}, and the recovery argument is the usual restricted isometry proof of stable $\ell^1$ reconstruction. See \cite[Chapters~6 and~12]{FoucartRauhut} for further background.

\section{Fourier ratio bounds on the periodic square}\label{sec:square}

Throughout this section, functions on $[0,1]^2$ are understood to extend periodically to $\mathbb R^2$, and we use the norm
$$ \|f\|_{C^1([0,1]^2)} = \max\left\{ \|f\|_{L^\infty}, \left\|\frac{\partial f}{\partial x_1}\right\|_{L^\infty}, \left\|\frac{\partial f}{\partial x_2}\right\|_{L^\infty} \right\}. $$
For $0<\gamma\leq1$ we also write
$$ [f]_{C^{0,\gamma}} = \sup_{u\neq v} \frac{|f(u)-f(v)|}{d_{\mathbb T^2}(u,v)^\gamma} \text{ and } \|f\|_{C^{0,\gamma}} = \|f\|_{L^\infty}+[f]_{C^{0,\gamma}}, $$
where $\mathbb T^2=\mathbb R^2/\mathbb Z^2$ is the torus, identified with the periodic square, and $d_{\mathbb T^2}$ is its distance. We begin with an estimate for Riemann sums.

\begin{lemma}[H\"older quadrature estimate]\label{lem:holder-riemann}
Let $0<\gamma\leq1$ and let $h\in C^{0,\gamma}([0,1]^2)$ be $1$-periodic in each variable. Then
$$ \left| \frac{1}{N^2}\sum_{x\in\mathbb Z_N^2}h\left(\frac{x}{N}\right) -\int_{[0,1]^2}h(u)du \right| \leq 2^{\frac{\gamma}{2}}N^{-\gamma}[h]_{C^{0,\gamma}}. $$
\end{lemma}

\begin{proof}
For each $x\in\mathbb Z_N^2$, let
$$ Q_x=\frac{x}{N}+[0,N^{-1})^2, $$
with addition understood periodically. These cells partition the torus, each of them has area $N^{-2}$, and every $u\in Q_x$ satisfies $d_{\mathbb T^2}(u,\frac{x}{N})\leq\frac{\sqrt2}{N}$. It follows that
$$ \begin{aligned}\left|\frac{1}{N^2}\sum_xh\left(\frac{x}{N}\right)-\int h(u)du\right| &= \left|\sum_x\int_{Q_x}\left(h\left(\frac{x}{N}\right)-h(u)\right)du\right|\\ &\leq \sum_x\int_{Q_x}[h]_{C^{0,\gamma}}d_{\mathbb T^2}\left(u,\frac{x}{N}\right)^\gamma du\\ &\leq 2^{\frac{\gamma}{2}}N^{-\gamma}[h]_{C^{0,\gamma}},\end{aligned} $$
which is the desired estimate.
\end{proof}

\subsection{H\"older and fractional Sobolev bounds}

The next estimate describes what remains of the discrete energy argument below the $C^1$ threshold. A H\"older increment across one grid cell has size $N^{-\gamma}$ rather than $N^{-1}$, and when $0<\gamma<1$ it is useful to match the size of the increment to the frequencies being estimated, which removes the logarithmic factor from the bound in that range.

\begin{proposition}\label{prop:holder-square-fr}
Let $N\geq2$ and $0<\gamma\leq1$, let $f$ be a nonzero function in $C^{0,\gamma}([0,1]^2)$ which is $1$-periodic in each variable, and define $g(x)=f(\frac{x}{N})$ for $x\in\mathbb Z_N^2$. Write
$$ A_\gamma(N)= \begin{cases} N^{1-\gamma},&0<\gamma<1,\\ \sqrt{1+\log N},&\gamma=1. \end{cases} $$
There exist constants $K_\gamma,C_\gamma>0$, depending only on $\gamma$, such that if
$$ N^\gamma\|f\|_{L^2([0,1]^2)}^2 \geq K_\gamma\|f\|_{C^{0,\gamma}([0,1]^2)}^2, $$
then
$$ FR(g)\leq\sqrt{2} \frac{\left|\int_{[0,1]^2}f(u)du\right|} {\|f\|_{L^2([0,1]^2)}} +C_\gamma \frac{\|f\|_{C^{0,\gamma}([0,1]^2)}} {\|f\|_{L^2([0,1]^2)}} (A_\gamma(N)+N^{-\gamma}). $$
\end{proposition}

\begin{proof}
We apply Lemma~\ref{lem:holder-riemann} to $h=|f|^2$. Since
$$ \left||f(u)|^2-|f(v)|^2\right| \leq(|f(u)|+|f(v)|)|f(u)-f(v)|, $$
we have
$$ [|f|^2]_{C^{0,\gamma}} \leq2\|f\|_{L^\infty}[f]_{C^{0,\gamma}} \leq2\|f\|_{C^{0,\gamma}}^2, $$
and therefore
$$ \frac{1}{N^2}\|g\|_{\ell^2}^2 \geq\|f\|_{L^2}^2 -2^{1+\frac{\gamma}{2}}N^{-\gamma}\|f\|_{C^{0,\gamma}}^2. $$
If $K_\gamma$ is chosen sufficiently large, it follows that
\begin{equation}\label{eq:holder-sampled-mass}
\|\widehat g\|_{\ell^2}=\|g\|_{\ell^2} \geq\frac{N}{\sqrt{2}}\|f\|_{L^2}.
\end{equation}
At the zero frequency we have
$$ \widehat g(0)=N\left(\frac{1}{N^2}\sum_x f\left(\frac{x}{N}\right)\right), $$
and Lemma~\ref{lem:holder-riemann} implies that
$$ |\widehat g(0)| \leq N\left|\int_{[0,1]^2}f(u)du\right| +2^{\frac{\gamma}{2}}N^{1-\gamma}[f]_{C^{0,\gamma}}. $$

We now turn to the nonzero frequencies, and we first suppose that $0<\gamma<1$. Recall that $|m_j|_*=\min\{m_j,N-m_j\}$, where $m_j$ is represented by an integer between $0$ and $N-1$, and put $|m|_\infty=\max\{|m_1|_*,|m_2|_*\}$. For a power of two $R$ with $1\leq R\leq\frac{N}{2}$, let
$$ \mathcal A_R= \{m\in\mathbb Z_N^2: R\leq|m|_\infty<2R\} $$
and choose the integer
$$ h_R=\max\left\{1,\left\lfloor\frac{N}{4R}\right\rfloor\right\}. $$
The sets $\mathcal A_R$ partition the nonzero frequencies, and each of them contains at most $CR^2$ points.

The point of this choice of $h_R$ is that the corresponding difference detects every frequency in $\mathcal A_R$. Indeed, if $m\in\mathcal A_R$, then $|m_j|_*\geq R$ for at least one $j$, and for that coordinate
$$ \frac{1}{8}\leq\frac{h_R|m_j|_*}{N}\leq\frac{1}{2}. $$
To see the lower bound, observe that $\lfloor\frac{N}{4R}\rfloor\geq\frac{N}{8R}$ when $R\leq\frac{N}{4}$, while $h_R=1$ when $R>\frac{N}{4}$, so the same lower bound holds. The upper bound follows from $|m_j|_*<2R$ in the first case and from $|m_j|_*\leq\frac{N}{2}$ in the second. Consequently,
$$ \sum_{j=1}^2 \left|e^{\frac{2\pi i h_Rm_j}{N}}-1\right|^2 \geq4\sin^2\left(\frac{\pi}{8}\right) \text{ for } m\in\mathcal A_R, $$
and Parseval's identity implies that
$$ \sum_{m\in\mathcal A_R}|\widehat g(m)|^2 \leq C\sum_{j=1}^2 \|g(\cdot+h_Re_j)-g\|_{\ell^2}^2. $$
By periodicity and H\"older continuity,
$$ \|g(\cdot+h_Re_j)-g\|_{\ell^2} \leq N\left(\frac{h_R}{N}\right)^\gamma [f]_{C^{0,\gamma}} \leq N R^{-\gamma}[f]_{C^{0,\gamma}}, $$
where we used $\frac{h_R}{N}\leq\frac{1}{R}$. It follows from Cauchy--Schwarz on $\mathcal A_R$ that
$$ \sum_{m\in\mathcal A_R}|\widehat g(m)| \leq C N R^{1-\gamma}[f]_{C^{0,\gamma}}, $$
and summing over the powers of two we conclude that
$$ \sum_{m\neq0}|\widehat g(m)| \leq C_\gamma N^{2-\gamma}[f]_{C^{0,\gamma}}. $$
The constant depends on $\gamma$ because the geometric sum has ratio $2^{1-\gamma}$, and it need not remain bounded as $\gamma$ approaches $1$.

When $\gamma=1$ we use instead the nearest neighbor differences $\Delta_jg(x)=g(x+e_j)-g(x)$, $j=1,2$, which satisfy $\|\Delta_jg\|_{\ell^2}\leq[f]_{C^{0,1}}$. By Proposition~\ref{prop:discrete-energy},
$$ \begin{aligned} \sum_{m\neq0}|\widehat g(m)| &\leq\left(\sum_{m\neq0}\frac{1}{\rho(m)^2}\right)^{\frac{1}{2}} (\|\Delta_1g\|_2^2+\|\Delta_2g\|_2^2)^{\frac{1}{2}}\\ &\leq C N\sqrt{1+\log N}[f]_{C^{0,1}}. \end{aligned} $$
In both cases the proposition follows by adding the zero frequency and dividing by the lower bound in \eqref{eq:holder-sampled-mass}.
\end{proof}

For $0<s<1$ and $1<p<\infty$ we use the periodic Sobolev--Slobodeckij norm
$$ \|f\|_{W^{s,p}}^p =\|f\|_{L^p}^p +\int_{[0,1]^2}\int_{[0,1]^2} \frac{|f(u)-f(v)|^p} {d_{\mathbb T^2}(u,v)^{2+sp}}dudv, $$
and we write $[f]_{W^{s,p}}$ for the periodic seminorm, the $p$-th root of the double integral above, so that $[f]_{W^{s,p}}\leq\|f\|_{W^{s,p}}$. The fractional Morrey inequality says that if $sp>2$ and $\alpha=s-\frac{2}{p}$, then
$$ \|f\|_{C^{0,\alpha}([0,1]^2)} \leq M_{s,p}\|f\|_{W^{s,p}([0,1]^2)}. $$
See, for example, \cite{AdamsFournier}. We shall also need a local version of this inequality, and for a ball $E\subset\mathbb R^2$ we write
$$ [f]_{W^{s,p}(E)}= \left(\int_E\int_E \frac{|f(u)-f(v)|^p}{|u-v|^{2+sp}}dudv\right)^{\frac{1}{p}}. $$

\begin{lemma}[Local fractional Morrey inequality]\label{lem:local-fractional-morrey}
Let $0<s<1$, $1<p<\infty$, $sp>2$, and put $\alpha=s-\frac{2}{p}$. There is a constant $M'_{s,p}$ such that, for every closed ball $B\subset\mathbb R^2$ of radius $r>0$ and every continuous $f$ on $B$ with $[f]_{W^{s,p}(B)}<\infty$,
\begin{equation}\label{eq:local-morrey}
|f(u)-f(v)|\leq M'_{s,p}r^\alpha[f]_{W^{s,p}(B)} \text{ for all } u,v\in B.
\end{equation}
\end{lemma}

\begin{proof}
We first take $B=B_1$ to be the closed unit ball centered at the origin, and we let $f_{B_1}$ denote the mean of $f$ over $B_1$. By Jensen's inequality and the fact that $|u-v|\leq2$,
$$ \begin{aligned} \|f-f_{B_1}\|_{L^p(B_1)}^p &\leq\frac{1}{|B_1|}\int_{B_1}\int_{B_1}|f(u)-f(v)|^pdvdu\\ &\leq \frac{2^{2+sp}}{\pi}[f]_{W^{s,p}(B_1)}^p. \end{aligned} $$
The fractional Morrey inequality on the unit ball \cite[Theorem~8.2]{DNPV}, applied to $f-f_{B_1}$, then shows that for $u,v\in B_1$
$$ \begin{aligned} |f(u)-f(v)| &\leq C_{s,p}(\|f-f_{B_1}\|_{L^p(B_1)}+[f]_{W^{s,p}(B_1)})|u-v|^\alpha\\ &\leq C'_{s,p}[f]_{W^{s,p}(B_1)}. \end{aligned} $$
For a ball $B$ with center $c$ and radius $r$, put $f_r(z)=f(c+rz)$. The substitution $u=c+rz$ and $v=c+rw$ shows that
$$ [f_r]_{W^{s,p}(B_1)}^p =r^{sp-2}[f]_{W^{s,p}(B)}^p =r^{\alpha p}[f]_{W^{s,p}(B)}^p, $$
and the estimate \eqref{eq:local-morrey} follows from the case of the unit ball applied to $f_r$.
\end{proof}

\begin{lemma}[Sobolev differences on the grid]\label{lem:sobolev-grid-differences}
Let $0<s<1$, $1<p<\infty$, and $sp>2$. Let $f$ be continuous and $1$-periodic with $[f]_{W^{s,p}}<\infty$, and define $g(x)=f(\frac{x}{N})$ for $x\in\mathbb Z_N^2$. For every integer $1\leq h\leq\frac{N}{4}$ and $j\in\{1,2\}$,
$$ \|g(\cdot+he_j)-g\|_{\ell^2} \leq 9^{\frac{1}{p}}M'_{s,p}N\left(\frac{h}{N}\right)^s[f]_{W^{s,p}}. $$
\end{lemma}

\begin{proof}
Since $s<1$ and $sp>2$, we have $p>2$. Represent $x\in\mathbb Z_N^2$ by a point of $\{0,\ldots,N-1\}^2$, and let $B_x\subset\mathbb R^2$ be the closed ball of radius $\frac{h}{N}$ centered at $\frac{x}{N}+\frac{h}{2N}e_j$. This ball contains $\frac{x}{N}$ and $\frac{x+he_j}{N}$, and its diameter is $\frac{2h}{N}\leq\frac{1}{2}$, so the quotient map $\mathbb R^2\to\mathbb T^2$ is injective on $B_x$ and $|u-v|=d_{\mathbb T^2}(u,v)$ for $u,v\in B_x$. Applying Lemma~\ref{lem:local-fractional-morrey} to the periodic extension of $f$ on $B_x$ we obtain
$$ |g(x+he_j)-g(x)|\leq M'_{s,p}\left(\frac{h}{N}\right)^\alpha[f]_{W^{s,p}(B_x)}, $$
where $\alpha=s-\frac{2}{p}$ as before.

A point $u\in\mathbb T^2$ lies in the image of $B_x$ only if $x$ lies within distance $h$ of $Nu-\frac{h}{2}e_j$ in $\mathbb R^2/N\mathbb Z^2$, and a disc of radius $h$ contains at most $(2h+1)^2\leq9h^2$ integer points. Consequently every pair $(u,v)\in\mathbb T^2\times\mathbb T^2$ is counted for at most $9h^2$ indices $x$, and
$$ \sum_{x\in\mathbb Z_N^2}[f]_{W^{s,p}(B_x)}^p \leq9h^2[f]_{W^{s,p}}^p. $$
Since $\alpha p=sp-2$, this implies that
$$ \sum_{x\in\mathbb Z_N^2}|g(x+he_j)-g(x)|^p \leq9(M'_{s,p})^p h^{sp}N^{2-sp}[f]_{W^{s,p}}^p. $$
By H\"older's inequality on the $N^2$ grid points, $\|v\|_{\ell^2}\leq N^{1-\frac{2}{p}}\|v\|_{\ell^p}$, and therefore
$$ \begin{aligned} \|g(\cdot+he_j)-g\|_{\ell^2} &\leq N^{1-\frac{2}{p}}\cdot9^{\frac{1}{p}}M'_{s,p}h^sN^{\frac{2}{p}-s}[f]_{W^{s,p}}\\ &=9^{\frac{1}{p}}M'_{s,p}N\left(\frac{h}{N}\right)^s[f]_{W^{s,p}}, \end{aligned} $$
which completes the proof.
\end{proof}

\begin{theorem}\label{thm:sobolev-square-fr}
Let $0<s<1$, $1<p<\infty$, and $sp>2$, and set $\alpha=s-\frac{2}{p}>0$. Let $f\in W^{s,p}([0,1]^2)$ be nonzero and $1$-periodic, using its continuous representative, and define $g(x)=f(\frac{x}{N})$. There exist constants $K_{s,p},C_{s,p}>0$ such that if
\begin{equation}\label{eq:sobolev-sampling-condition}
N^\alpha\|f\|_{L^2([0,1]^2)}^2 \geq K_{s,p}\|f\|_{W^{s,p}([0,1]^2)}^2,
\end{equation}
then
\begin{equation}\label{eq:sobolev-sharp-fr}
FR(g)\leq\sqrt2 \frac{\left|\int_{[0,1]^2}f(u)du\right|}{\|f\|_{L^2([0,1]^2)}} +C_{s,p}\frac{\|f\|_{W^{s,p}([0,1]^2)}}{\|f\|_{L^2([0,1]^2)}}N^{1-s}.
\end{equation}
\end{theorem}

\begin{proof}
Since $p>2$ and $[0,1]^2$ has measure one, $\|f\|_{L^2}\leq\|f\|_{L^p}\leq\|f\|_{W^{s,p}}$. Let $K_\alpha$ be the constant of Proposition~\ref{prop:holder-square-fr} for $\gamma=\alpha$, let $M_{s,p}$ be the constant in the global fractional Morrey inequality above, and choose $K_{s,p}\geq\max\{4,K_\alpha M_{s,p}^2\}$. Then the condition \eqref{eq:sobolev-sampling-condition} implies that $N^\alpha\geq K_{s,p}\geq4$, hence $N\geq4$ because $\alpha<1$. The global Morrey inequality and the condition \eqref{eq:sobolev-sampling-condition} also imply the size condition of Proposition~\ref{prop:holder-square-fr} for $\gamma=\alpha$, so the first part of the proof of that proposition shows that
\begin{equation}\label{eq:sobolev-mass-and-zero}
\begin{gathered}
\|\widehat g\|_{\ell^2}=\|g\|_{\ell^2}\geq\frac{N}{\sqrt2}\|f\|_{L^2},\\
|\widehat g(0)|\leq N\left|\int_{[0,1]^2}f(u)du\right| +2^{\frac{\alpha}{2}}M_{s,p}N^{1-\alpha}\|f\|_{W^{s,p}}.
\end{gathered}
\end{equation}
For the nonzero frequencies we follow the dyadic argument in the proof of Proposition~\ref{prop:holder-square-fr}, with $\mathcal A_R$ and $h_R$ as defined there. Since $N\geq4$, we have $1\leq h_R\leq\frac{N}{4}$, so Lemma~\ref{lem:sobolev-grid-differences} applies with $h=h_R$, and together with the lower bound for the Fourier multiplier from that proof and the inequality $\frac{h_R}{N}\leq\frac{1}{R}$ it shows that
$$ \sum_{m\in\mathcal A_R}|\widehat g(m)|^2 \leq C\sum_{j=1}^2\|g(\cdot+h_Re_j)-g\|_{\ell^2}^2 \leq C_{s,p}N^2R^{-2s}[f]_{W^{s,p}}^2. $$
Since $\mathcal A_R$ has at most $CR^2$ points, it follows from Cauchy--Schwarz that
$$ \sum_{m\in\mathcal A_R}|\widehat g(m)| \leq C_{s,p}NR^{1-s}[f]_{W^{s,p}}. $$
Summing over the powers of two $R\leq\frac{N}{2}$, which is a geometric sum with ratio $2^{1-s}>1$, we obtain
$$ \sum_{m\ne0}|\widehat g(m)| \leq\frac{C_{s,p}}{1-2^{s-1}}N^{2-s}[f]_{W^{s,p}}. $$
The estimate \eqref{eq:sobolev-sharp-fr} now follows by adding the zero frequency from \eqref{eq:sobolev-mass-and-zero}, dividing by the lower bound for $\|\widehat g\|_{\ell^2}$, and using $[f]_{W^{s,p}}\leq\|f\|_{W^{s,p}}$ and $N^{-\alpha}\leq1\leq N^{1-s}$.
\end{proof}

Proposition~\ref{prop:holder-square-fr} applies in the same way to any Sobolev class of higher order which the Sobolev--Morrey embedding places in a H\"older class with a positive exponent. Theorem~\ref{thm:sobolev-square-fr} deals with the fractional range, where the condition $sp>2$ guarantees a continuous representative and the local difference argument is available.

\begin{remark}\label{rem:holder-threshold}
The trivial estimate on the group $\mathbb Z_N^2$, which has $N^2$ points, is $FR(g)\leq N$. Since $N^{1-\gamma}=o(N)$ for every fixed $0<\gamma<1$, Proposition~\ref{prop:holder-square-fr} improves on this estimate as soon as the function has any positive H\"older exponent, and Theorem~\ref{thm:sobolev-square-fr} does so for every $0<s<1$ in the range $sp>2$. We shall see in Remark~\ref{rem:exponents-sharp} that neither polynomial exponent can be lowered in estimates that are uniform over the corresponding regularity classes. At $\gamma=1$ the proposition yields the critical $\sqrt{1+\log N}$ bound, which therefore holds for every $C^1$ function.
\end{remark}

\begin{remark}[Sharpness of the exponents]\label{rem:exponents-sharp}
Let $0<\theta<1$, let $N\geq3$ be odd, and put
$$ \psi(x)=\cos\left(\frac{2\pi(x_1^2+x_2^2)}{N}\right) \text{ for } x\in\mathbb Z_N^2. $$
Let $\mathsf H(\tau)=\max\{0,1-|\tau|\}$, extend $\psi$ $N$-periodically to $\mathbb Z^2$, and let
$$ \Phi(u)=\sum_{x\in\mathbb Z^2}\psi(x)\mathsf H(Nu_1-x_1)\mathsf H(Nu_2-x_2) \text{ and } f_N=1+N^{-\theta}\Phi, $$
so that $\Phi$ is the bilinear interpolant of $\psi$ and $f_N$ is the corresponding perturbation of the constant function. Then $f_N$ is real, $1$-periodic, and Lipschitz, and it satisfies $\|f_N\|_{C^{0,\theta}}\leq5$, $\|f_N\|_{W^{\theta,p}}\leq C_{\theta,p}$ for every $1<p<\infty$, and $\frac{2}{3}\leq\|f_N\|_{L^2}\leq2$, while $g_N(x)=f_N(\frac{x}{N})$ satisfies
\begin{equation}\label{eq:chirp-sharp-fr}
FR(g_N)\geq\frac{N^{1-\theta}}{4}-\frac{1}{2}.
\end{equation}
For large $N$ the size conditions of Proposition~\ref{prop:holder-square-fr} and Theorem~\ref{thm:sobolev-square-fr} hold, so taking $\theta=\gamma$ and $\theta=s$ we see that the exponents $1-\gamma$ and $1-s$ in these statements are sharp for bounds that are uniform over such families of functions.

To verify the estimate \eqref{eq:chirp-sharp-fr}, recall that
$$ \left|\sum_{y\in\mathbb Z_N}e^{\frac{2\pi i(ay^2+by)}{N}}\right|=\sqrt N $$
when $N$ is odd and $\gcd(a,N)=1$. It follows that the function $\chi(x)=e^{\frac{2\pi i(x_1^2+x_2^2)}{N}}$ satisfies $|\widehat\chi(m)|=1$ for every $m$ and
$$ \left|\sum_{x\in\mathbb Z_N^2}\chi(x)^2\right|=N. $$
Since $\psi=\frac{\chi+\overline\chi}{2}$, we have $|\widehat\psi(m)|\leq1$ for every $m$ and
$$ \|\psi\|_{\ell^2}^2 =\frac{N^2}{2}+\frac{1}{2}\operatorname{Re}\sum_{x\in\mathbb Z_N^2}\chi(x)^2 \geq\frac{N^2-N}{2}, $$
and therefore
$$ \sum_{m\ne0}|\widehat\psi(m)| \geq\sum_{m\ne0}|\widehat\psi(m)|^2 \geq\frac{N^2-N}{2}-1. $$
Since $g_N=1+N^{-\theta}\psi$, we have $\widehat g_N(m)=N^{-\theta}\widehat\psi(m)$ for $m\ne0$ and $\|g_N\|_{\ell^2}\leq(1+N^{-\theta})N\leq2N$, which proves \eqref{eq:chirp-sharp-fr}.

On the continuous side, $\Phi(u)$ is a convex combination of values of $\psi$, so $|\Phi|\leq1$. On each grid cell $\Phi$ is bilinear with $|\partial_j\Phi|\leq2N$, hence
$$ |\Phi(u)-\Phi(v)|\leq\min\{2,2\sqrt2Nd_{\mathbb T^2}(u,v)\}, $$
and consequently $[\Phi]_{C^{0,\theta}}\leq2\sqrt2N^\theta$. Splitting the double integral at $d_{\mathbb T^2}(u,v)=N^{-1}$ we see in the same way that $[\Phi]_{W^{\theta,p}}\leq C_{\theta,p}N^\theta$, and these estimates imply the stated bounds for $\|f_N\|_{C^{0,\theta}}$ and $\|f_N\|_{W^{\theta,p}}$. Finally, the integral
$$ \int_{[0,1]^2}\Phi(u)du =N^{-2}\sum_{x\in\mathbb Z_N^2}\psi(x) $$
has absolute value at most $N^{-1}$, so $\|f_N\|_{L^2}\geq\left|\int f_N\right|\geq1-N^{-1-\theta}\geq\frac{2}{3}$ and $\|f_N\|_{L^2}\leq\|f_N\|_{L^\infty}\leq2$.
\end{remark}

\subsection{A translated grid}

The H\"older estimate on a fixed grid requires a quantitative hypothesis which keeps the sampled $L^2$ norm from becoming too small. For $C^1$ functions a translation of the grid removes this issue, and averaging over the translation parameter allows us to replace the $C^1$ norm in the main term by the $L^2$ norm of the gradient.

\begin{theorem}\label{thm:translated-grid}
Let $f$ be a nonzero function in $C^1([0,1]^2)$ which is $1$-periodic in each variable, let $N\geq2$, and for $t=(t_1,t_2)\in[0,1)^2$ let
$$ g_t(x_1,x_2) = f\left(\frac{x_1+t_1}{N},\frac{x_2+t_2}{N}\right). $$
There exists $t\in[0,1)^2$ such that
$$ FR(g_t) \leq 1+ \frac{1}{\sqrt{2}}\sqrt{1+\log N} \frac{\|\nabla f\|_{L^2([0,1]^2)}}{\|f\|_{L^2([0,1]^2)}}. $$
There is also a translation $t'\in[0,1)^2$, not necessarily equal to $t$, such that
$$ FR(g_{t'})^2 \leq 2\frac{\left|\int_{[0,1]^2}f(u)du\right|^2}{\|f\|_{L^2([0,1]^2)}^2} + \left(1+\log N+\frac{1}{2\pi^2N^2}\right) \frac{\|\nabla f\|_{L^2([0,1]^2)}^2}{\|f\|_{L^2([0,1]^2)}^2}. $$
If $f\geq 0$ and $\mu=\int_{[0,1]^2}f(u)du>0$, then
$$ FR(g_t) \leq 1+ \frac{1}{\sqrt{2}}\sqrt{1+\log N} \frac{\|\nabla f\|_{L^2([0,1]^2)}}{\mu}. $$
\end{theorem}

\begin{proof}
For $t\in[0,1)^2$ we let $D(t)=\|g_t\|_2^2$, $E(t)=\|\Delta_1g_t\|_2^2+\|\Delta_2g_t\|_2^2$, and $Z(t)=|\widehat{g_t}(0)|$, and we also write
$$ W_N=\sum_{m\neq0}\frac{1}{\rho(m)^2}. $$
By Proposition~\ref{prop:discrete-energy}, for every $t$ with $D(t)>0$,
$$ FR(g_t) \leq \frac{Z(t)+\sqrt{W_N}\sqrt{E(t)}}{\sqrt{D(t)}}, $$
and the elementary inequality $(a+b)^2\leq2a^2+2b^2$ shows that
$$ FR(g_t)^2 \leq \frac{2Z(t)^2+2W_NE(t)}{D(t)}. $$
The idea is to average the numerator and the denominator of this last expression over $t$. For each fixed $x$, the change of variables $u=\frac{x+t}{N}$ maps the translation square onto the grid cell indexed by $x$ and has Jacobian $dt=N^2du$, so summing over the cells we obtain
$$ \int_{[0,1)^2}D(t)dt = N^2\|f\|_{L^2([0,1]^2)}^2. $$
In the same way, for each $j$,
$$ \int_{[0,1)^2}\|\Delta_jg_t\|_2^2dt = N^2\int_{[0,1]^2} \left|f\left(u+\frac{e_j}{N}\right)-f(u)\right|^2du, $$
and by the fundamental theorem of calculus and Cauchy--Schwarz,
$$ \left|f\left(u+\frac{e_j}{N}\right)-f(u)\right|^2 = \left|\int_0^{\frac{1}{N}} \frac{\partial f}{\partial x_j}(u+re_j)dr\right|^2 \leq \frac{1}{N}\int_0^{\frac{1}{N}} \left|\frac{\partial f}{\partial x_j}(u+re_j)\right|^2dr. $$
Integrating in $u$ and using the translation invariance of the torus we conclude that
$$ \int_{[0,1)^2}E(t)dt \leq \|\nabla f\|_{L^2([0,1]^2)}^2. $$
These two averaged estimates imply that there is a translation $t_0$ with $D(t_0)>0$ and
$$ \frac{E(t_0)}{D(t_0)} \leq \frac{\|\nabla f\|_{L^2([0,1]^2)}^2}{N^2\|f\|_{L^2([0,1]^2)}^2}, $$
because if the reverse strict inequality held whenever $D(t)>0$, then multiplication by $D(t)$ and integration would contradict them. Since a single Fourier coefficient is at most the full $\ell^2$ norm, we have $Z(t_0)\leq\sqrt{D(t_0)}$, and the estimate for $W_N$ from Proposition~\ref{prop:discrete-energy} shows that
$$ FR(g_{t_0}) \leq 1+\sqrt{W_N}\sqrt{\frac{E(t_0)}{D(t_0)}} \leq 1+ \frac{1}{\sqrt{2}}\sqrt{1+\log N} \frac{\|\nabla f\|_{L^2([0,1]^2)}}{\|f\|_{L^2([0,1]^2)}}. $$
This is the first assertion, and the bound for nonnegative $f$ follows because $\mu\leq\|f\|_{L^2([0,1]^2)}$ when $f\geq0$.

In order to obtain a version which retains the continuous mean, we also average the zero frequency. Let
$$ c(k)=\int_{[0,1]^2}f(u)e^{-2\pi i k\cdot u}du \text{ for } k\in\mathbb Z^2 $$
be the continuous Fourier coefficients of $f$. We claim that the function
$$ t\mapsto\widehat{g_t}(0) = \frac{1}{N}\sum_{x\in\mathbb Z_N^2}f\left(\frac{x+t}{N}\right) $$
has continuous Fourier coefficient $Nc(N\ell)$ at the frequency $\ell\in\mathbb Z^2$. To verify this directly, let $Q_x$ be the grid cell indexed by $x$, as in the proof of Lemma~\ref{lem:holder-riemann}. Then
$$ \begin{aligned} \int_{[0,1)^2}\widehat{g_t}(0)e^{-2\pi i\ell\cdot t}dt &= \frac{1}{N}\sum_{x\in\mathbb Z_N^2} N^2\int_{Q_x}f(u)e^{-2\pi i\ell\cdot(Nu-x)}du\\ &= N\int_{[0,1]^2}f(u)e^{-2\pi iN\ell\cdot u}du = Nc(N\ell), \end{aligned} $$
where we used $e^{2\pi i\ell\cdot x}=1$ in the second step. By Parseval's identity in the variable $t$,
$$ \int_{[0,1)^2}Z(t)^2dt = N^2\sum_{\ell\in\mathbb Z^2}|c(N\ell)|^2 \leq N^2\left|\int_{[0,1]^2}f(u)du\right|^2 + \frac{1}{4\pi^2}\|\nabla f\|_{L^2([0,1]^2)}^2, $$
where the inequality holds because $|N\ell|\geq N$ when $\ell\neq0$ and
$$ 4\pi^2\sum_{k\in\mathbb Z^2}|k|^2|c(k)|^2 = \|\nabla f\|_{L^2([0,1]^2)}^2. $$
Recall from Proposition~\ref{prop:discrete-energy} that $W_N\leq\frac{N^2}{2}(1+\log N)$. Since the integral of $D(t)$ is positive, there is a translation $t'$ with $D(t')>0$ such that
$$ \frac{2Z(t')^2+2W_NE(t')}{D(t')} \leq \frac{ \int_{[0,1)^2}(2Z(t)^2+2W_NE(t))dt }{ \int_{[0,1)^2}D(t)dt }, $$
since otherwise the strict reverse inequality would remain true after multiplication by $D(t)$ and integration. The squared bound in the theorem follows by substituting the three averaged estimates.
\end{proof}

\begin{corollary}[A $C^1$ estimate on a translated grid]\label{cor:translated-grid-c1}
Under the hypotheses of Theorem~\ref{thm:translated-grid}, there exists a translation $\tau\in[0,1)^2$ such that
$$ FR(g_\tau) \leq \frac{\left|\int_{[0,1]^2}f(u)du\right|}{\|f\|_{L^2([0,1]^2)}} + \frac{\|f\|_{C^1([0,1]^2)}}{\|f\|_{L^2([0,1]^2)}}\sqrt{1+\log N} + \frac{2}{N} \frac{\|f\|_{C^1([0,1]^2)}}{\|f\|_{L^2([0,1]^2)}}. $$
If $f\geq0$ and $\mu=\int_{[0,1]^2}f(u)du>0$, then
$$ FR(g_\tau) \leq 1+ \frac{\|f\|_{C^1([0,1]^2)}}{\mu}\sqrt{1+\log N} + \frac{2}{N} \frac{\|f\|_{C^1([0,1]^2)}}{\mu}. $$
\end{corollary}

\begin{proof}
The averaged identity
$$ \int_{[0,1)^2}\|g_t\|_2^2dt = N^2\|f\|_{L^2([0,1]^2)}^2 $$
shows that some translation $\tau$ satisfies $\|g_\tau\|_2\geq N\|f\|_{L^2([0,1]^2)}$. On the other hand, for every translation, Lemma~\ref{lem:holder-riemann} with $\gamma=1$ and the estimate for the discrete differences coming from the fundamental theorem of calculus show that
$$ |\widehat{g_\tau}(0)| \leq N\left|\int_{[0,1]^2}f(u)du\right| +2\|f\|_{C^1([0,1]^2)} $$
and
$$ \sum_{m\neq0}|\widehat{g_\tau}(m)| \leq N\|f\|_{C^1([0,1]^2)}\sqrt{1+\log N}. $$
The first assertion follows by adding these estimates and dividing by the lower bound for $\|g_\tau\|_2=\|\widehat{g_\tau}\|_2$, and the second assertion follows from the first because $\mu\leq\|f\|_{L^2([0,1]^2)}$ when $f\geq0$.
\end{proof}

\subsection{Two derivatives}

Two derivatives remove the logarithmic factor in the estimate on the square, and the proof is the same discrete energy argument with second differences in place of first differences. We use the norm
$$ \|f\|_{C^2([0,1]^2)} =\max_{\beta_1+\beta_2\leq2} \left\|\frac{\partial^{\beta_1+\beta_2}f} {\partial x_1^{\beta_1}\partial x_2^{\beta_2}}\right\|_{L^\infty([0,1]^2)}, $$
where $\beta_1$ and $\beta_2$ are nonnegative integers and the derivative of order zero is $f$ itself, and we write $\Delta f=\frac{\partial^2f}{\partial x_1^2}+\frac{\partial^2f}{\partial x_2^2}$ for the continuous Laplacian.

\begin{proposition}\label{prop:c2-square-fr}
Let $f$ be a nonzero $C^2$ function on $[0,1]^2$ which is $1$-periodic in each variable, let $N\geq2$, and let $g(x)=f(\frac{x}{N})$ for $x\in\mathbb Z_N^2$. There are absolute constants $K,C>0$ such that
$$ N\|f\|_{L^2([0,1]^2)}^2 \geq K\|f\|_{C^2([0,1]^2)}^2 $$
implies
$$ FR(g)\leq1+C \frac{\|f\|_{C^2([0,1]^2)}} {\|f\|_{L^2([0,1]^2)}}. $$
For every $N\geq2$, without this size condition, there exists a translation $t\in[0,1)^2$ such that $g_t(x)=f(\frac{x+t}{N})$ is nonzero and
$$ FR(g_t)\leq1+C \frac{\|\Delta f\|_{L^2([0,1]^2)}} {\|f\|_{L^2([0,1]^2)}}. $$
\end{proposition}

\begin{proof}
Define the discrete Laplacian with nonnegative Fourier multiplier by
$$ \mathcal L_Ng(x)=\sum_{j=1}^2 (2g(x)-g(x+e_j)-g(x-e_j)). $$
Its Fourier transform is $\rho(m)^2\widehat g(m)$, so Proposition~\ref{prop:discrete-energy} with $s=2$ implies that
\begin{equation}\label{eq:c2-discrete-energy}
\sum_{m\neq0}|\widehat g(m)| \leq\left(\sum_{m\neq0}\frac{1}{\rho(m)^4}\right)^{\frac{1}{2}} \|\mathcal L_Ng\|_2 \leq C N^2\|\mathcal L_Ng\|_2.
\end{equation}
Applying the fundamental theorem of calculus twice we see that
$$ \left|2f(u)-f\left(u+\frac{e_j}{N}\right) -f\left(u-\frac{e_j}{N}\right)\right| \leq\frac{1}{N^2} \left\|\frac{\partial^2 f}{\partial x_j^2}\right\|_{L^\infty}, $$
and since there are $N^2$ grid points, $\|\mathcal L_Ng\|_2\leq\frac{C}{N}\|f\|_{C^2}$. The size condition and the argument for the sampled mass in the proof of Proposition~\ref{prop:holder-square-fr} show that $\|g\|_2\geq\frac{N}{\sqrt{2}}\|f\|_{L^2}$. The assertion on the fixed grid follows by dividing the estimate \eqref{eq:c2-discrete-energy} by this lower bound and using $|\widehat g(0)|\leq\|g\|_2$.

For the translated assertion, let $D(t)=\|g_t\|_2^2$ and $E_2(t)=\|\mathcal L_Ng_t\|_2^2$. The same change of variables as in the proof of Theorem~\ref{thm:translated-grid} shows that
$$ \int_{[0,1)^2}D(t)dt=N^2\|f\|_{L^2}^2 $$
and
$$ \int_{[0,1)^2}E_2(t)dt =N^2\left\|\sum_{j=1}^2 \left(2f(\cdot)-f\left(\cdot+\frac{e_j}{N}\right) -f\left(\cdot-\frac{e_j}{N}\right)\right) \right\|_{L^2}^2. $$
In order to estimate the last expression, let
$$ c(k)=\int_{[0,1]^2}f(u)e^{-2\pi i k\cdot u}du \text{ for } k\in\mathbb Z^2 $$
as before. The Fourier multiplier of the sum of continuous differences is
$$ 4\sum_{j=1}^2\sin^2\left(\frac{\pi k_j}{N}\right) \leq\frac{4\pi^2}{N^2}|k|^2, $$
and it follows from Parseval's identity that
$$ \int_{[0,1)^2}E_2(t)dt \leq\frac{1}{N^2}\|\Delta f\|_{L^2}^2. $$
Averaging $E_2(t)$ and $D(t)$ together as in the proof of Theorem~\ref{thm:translated-grid}, we obtain a translation $t$ with $D(t)>0$ and
$$ \frac{E_2(t)}{D(t)} \leq\frac{1}{N^4} \frac{\|\Delta f\|_{L^2}^2}{\|f\|_{L^2}^2}, $$
and the proof is completed by applying \eqref{eq:c2-discrete-energy} to $g_t$ and using once again that $|\widehat{g_t}(0)|\leq\|g_t\|_2$.
\end{proof}

\section{Fourier ratio bounds on compact manifolds}\label{sec:manifolds}

Let $(M,g)$ be a connected compact $d$-dimensional Riemannian manifold without boundary, and write $|M|$ for its Riemannian volume. We fix an orthonormal basis $\{e_j\}_{j=0}^\infty$ of $L^2(M)$ consisting of eigenfunctions of the Laplace--Beltrami operator, so that
$$ -\Delta_ge_j=\lambda_je_j \text{ and } 0=\lambda_0<\lambda_1\leq\lambda_2\leq\cdots, $$
and we take $e_0=\frac{1}{\sqrt{|M|}}$. For $\Lambda\geq 1$, let
$$ V_\Lambda=\operatorname{span}\{e_j:\lambda_j\leq \Lambda\} \text{ and } D_\Lambda=\#\{j:\lambda_j\leq \Lambda\}, $$
and for a nonzero $f\in V_\Lambda$ define
$$ \widehat f(j)=\int_Mf(x)\overline{e_j(x)}dV_g(x) \text{ and } FR_\Lambda(f)=\frac{\sum_{\lambda_j\leq \Lambda}|\widehat f(j)|}{\left(\sum_{\lambda_j\leq \Lambda}|\widehat f(j)|^2\right)^{\frac{1}{2}}}. $$
The spectral decomposition also defines the fractional powers of the Laplace--Beltrami operator. If $s>0$ and $f\in V_\Lambda$, then
$$ (-\Delta_g)^{\frac{s}{2}}f = \sum_{0<\lambda_j\leq \Lambda} \lambda_j^{\frac{s}{2}}\widehat f(j)e_j, $$
where the constant mode is omitted because its eigenvalue is zero, and by Parseval's identity
$$ \|(-\Delta_g)^{\frac{s}{2}}f\|_{L^2(M)}^2 = \sum_{0<\lambda_j\leq \Lambda}\lambda_j^s|\widehat f(j)|^2. $$
When $s=1$, integration by parts on a compact manifold without boundary leads to the Dirichlet energy identity
$$ \|(-\Delta_g)^{\frac{1}{2}}f\|_{L^2(M)}^2 = \langle-\Delta_gf,f\rangle_{L^2(M)} = \|\nabla_gf\|_{L^2(M)}^2, $$
so the following proposition is the spectral counterpart of the discrete gradient calculation on $\mathbb Z_N^2$.

\begin{proposition}\label{prop:manifold-fr}
Let $f\in V_\Lambda$ be nonzero. For every $s>0$,
$$ FR_\Lambda(f) \leq \frac{1}{\sqrt{|M|}} \frac{\left|\int_Mf(x)dV_g(x)\right|}{\|f\|_{L^2(M)}} + \left(\sum_{0<\lambda_j\leq \Lambda}\frac{1}{\lambda_j^s}\right)^{\frac{1}{2}} \frac{\|(-\Delta_g)^{\frac{s}{2}}f\|_{L^2(M)}}{\|f\|_{L^2(M)}}. $$
Consequently,
$$ FR_\Lambda(f) \leq \frac{1}{\sqrt{|M|}} \frac{\left|\int_Mf(x)dV_g(x)\right|}{\|f\|_{L^2(M)}} + C_{M,s} \frac{\|(-\Delta_g)^{\frac{s}{2}}f\|_{L^2(M)}}{\|f\|_{L^2(M)}} \begin{cases} 1, & s>\frac{d}{2},\\ (1+\log \Lambda)^{\frac{1}{2}}, & s=\frac{d}{2},\\ \Lambda^{\frac{d}{4}-\frac{s}{2}}, & 0<s<\frac{d}{2}. \end{cases} $$
\end{proposition}

\begin{proof}
By Parseval's identity,
$$ \|f\|_{L^2(M)}^2 = \sum_{\lambda_j\leq \Lambda}|\widehat f(j)|^2, $$
and the zero mode satisfies
$$ |\widehat f(0)| = \frac{1}{\sqrt{|M|}} \left|\int_Mf(x)dV_g(x)\right|. $$
For the remaining modes we use weighted Cauchy--Schwarz to see that
$$ \begin{aligned} \sum_{0<\lambda_j\leq \Lambda}|\widehat f(j)| &\leq \left(\sum_{0<\lambda_j\leq \Lambda}\frac{1}{\lambda_j^s}\right)^{\frac{1}{2}} \left(\sum_{0<\lambda_j\leq \Lambda}\lambda_j^s|\widehat f(j)|^2\right)^{\frac{1}{2}}\\ &= \left(\sum_{0<\lambda_j\leq \Lambda}\frac{1}{\lambda_j^s}\right)^{\frac{1}{2}} \|(-\Delta_g)^{\frac{s}{2}}f\|_{L^2(M)}, \end{aligned} $$
which proves the first estimate.

Let us briefly recall how Weyl's law enters. Its full asymptotic form is
$$ \#\{j:\lambda_j\leq R\} = c_MR^{\frac{d}{2}}+o(R^{\frac{d}{2}}) \text{ as } R\to\infty, $$
but only the resulting upper bound
$$ \#\{j:\lambda_j\leq R\}\leq C_M(1+R)^{\frac{d}{2}} $$
is needed here. The finitely many eigenvalues below $1$ can be absorbed into the constant. On the dyadic shell $2^k\leq\lambda_j<2^{k+1}$ there are at most $C_M2^{\frac{kd}{2}}$ eigenvalues, while each weight $\frac{1}{\lambda_j^s}$ is at most $\frac{1}{2^{ks}}$, so the contribution of this shell is bounded by $C_M2^{k(\frac{d}{2}-s)}$. These terms form a convergent geometric series if $s>\frac{d}{2}$, each shell contributes a bounded amount and there are $O(1+\log \Lambda)$ shells if $s=\frac{d}{2}$, and the final shell dominates and contributes $\Lambda^{\frac{d}{2}-s}$ if $s<\frac{d}{2}$, so the three cases in the statement follow upon taking square roots. See \cite{Sogge} for the spectral estimates used here.
\end{proof}

\begin{corollary}\label{cor:surface-fr}
If $M$ is two-dimensional and $f\in V_\Lambda$ is nonzero, then
$$ FR_\Lambda(f) \leq \frac{1}{\sqrt{|M|}} \frac{\left|\int_Mf(x)dV_g(x)\right|}{\|f\|_{L^2(M)}} + C_M(1+\log \Lambda)^{\frac{1}{2}} \frac{\|\nabla_gf\|_{L^2(M)}}{\|f\|_{L^2(M)}}. $$
\end{corollary}

\subsection{Sharpness on compact surfaces}

It is natural to ask whether the square-root logarithm in Corollary~\ref{cor:surface-fr} is an artifact of the proof. The next result shows that it is not, and that it is forced by the distribution of the eigenvalues.

\begin{theorem}\label{thm:surface-sharpness}
Let $M$ be a connected compact two-dimensional Riemannian manifold without boundary, and fix an orthonormal eigenfunction basis as above. For every sufficiently large $\Lambda$, define
$$ f_\Lambda = \sum_{0<\lambda_j\leq \Lambda}\frac{1}{\lambda_j}e_j. $$
There are positive constants $c_M$ and $C_M$, independent of $\Lambda$, such that
$$ c_M\log \Lambda \leq FR_\Lambda(f_\Lambda) \leq C_M\log \Lambda $$
and
$$ c_M\sqrt{\log \Lambda} \leq \frac{\|\nabla_gf_\Lambda\|_{L^2(M)}}{\|f_\Lambda\|_{L^2(M)}} \leq C_M\sqrt{\log \Lambda}. $$
Consequently,
$$ \frac{FR_\Lambda(f_\Lambda)\|f_\Lambda\|_{L^2(M)}}{\|\nabla_gf_\Lambda\|_{L^2(M)}} \asymp \sqrt{\log \Lambda}. $$
In particular, the factor $\sqrt{\log \Lambda}$ in the critical energy estimate cannot be replaced by a function that is $o(\sqrt{\log \Lambda})$.
\end{theorem}

\begin{proof}
Let
$$ S_1(\Lambda)=\sum_{0<\lambda_j\leq \Lambda}\frac{1}{\lambda_j} \text{ and } S_2(\Lambda)=\sum_{0<\lambda_j\leq \Lambda}\frac{1}{\lambda_j^2}. $$
The coefficient of $e_j$ in $f_\Lambda$ is $\frac{1}{\lambda_j}$ for every positive eigenvalue below $\Lambda$, so $\|\widehat{f_\Lambda}\|_1=S_1(\Lambda)$ and $\|f_\Lambda\|_{L^2(M)}^2=S_2(\Lambda)$, and by the Dirichlet energy identity
$$ \|\nabla_gf_\Lambda\|_{L^2(M)}^2 = \sum_{0<\lambda_j\leq \Lambda} \lambda_j\left|\frac{1}{\lambda_j}\right|^2 = S_1(\Lambda). $$
We now estimate the two sums. Weyl's law in dimension two implies that $\lambda_j\asymp j$ for the positive eigenvalues, once the index is sufficiently large, and that the number of eigenvalues below $\Lambda$ is comparable to $\Lambda$. Therefore
$$ S_1(\Lambda) \asymp \sum_{1\leq j\leq c\Lambda}\frac{1}{j} \asymp \log \Lambda. $$
On the other hand, the convergence of $\sum_{j=1}^{\infty}\frac{1}{j^2}$ shows that $S_2(\Lambda)$ is bounded above independently of $\Lambda$, and it is bounded below by a positive constant because it contains the contribution of the first positive eigenvalue, so $S_2(\Lambda)\asymp1$. We have shown that
$$ FR_\Lambda(f_\Lambda) = \frac{S_1(\Lambda)}{\sqrt{S_2(\Lambda)}} \asymp \log \Lambda $$
and
$$ \frac{\|\nabla_gf_\Lambda\|_{L^2(M)}}{\|f_\Lambda\|_{L^2(M)}} = \frac{\sqrt{S_1(\Lambda)}}{\sqrt{S_2(\Lambda)}} \asymp \sqrt{\log \Lambda}, $$
and dividing the first estimate by the second completes the proof.
\end{proof}

\subsection{The sphere}

On $S^2$ the eigenvalues are $\lambda_\ell=\ell(\ell+1)$ with multiplicity $2\ell+1$, so the critical estimate can be written with an explicit constant. Here $\sigma$ denotes the surface measure on $S^2$, and $\|f\|_{C^1(S^2)}$ is the maximum of the $L^\infty$ norms of $f$ and its spherical gradient.

\begin{corollary}\label{cor:sphere-fr}
Let $L\geq1$ be an integer, let $\mathcal H_L=V_{L(L+1)}$ be the space of spherical harmonics of degree at most $L$ on $S^2$, and use the standard orthonormal basis $\{Y_\ell^m\}$. For every nonzero $f\in\mathcal H_L$,
$$ FR_{L(L+1)}(f) \leq \frac{1}{\sqrt{4\pi}} \frac{\left|\int_{S^2}f(\omega)d\sigma(\omega)\right|}{\|f\|_{L^2(S^2)}} + \sqrt{2(1+\log L)} \frac{\|\nabla_{S^2}f\|_{L^2(S^2)}}{\|f\|_{L^2(S^2)}}. $$
In particular,
$$ FR_{L(L+1)}(f) \leq \frac{1}{\sqrt{4\pi}} \frac{\left|\int_{S^2}f(\omega)d\sigma(\omega)\right|}{\|f\|_{L^2(S^2)}} + \sqrt{8\pi(1+\log L)} \frac{\|f\|_{C^1(S^2)}}{\|f\|_{L^2(S^2)}}. $$
If $f\geq 0$ and $\mu=\int_{S^2}f d\sigma>0$, then
$$ FR_{L(L+1)}(f) \leq 1+4\sqrt{2}\pi \frac{\|f\|_{C^1(S^2)}}{\mu} \sqrt{1+\log L}. $$
\end{corollary}

\begin{proof}
Write
$$ f=\sum_{\ell=0}^L\sum_{m=-\ell}^{\ell} \widehat f(\ell,m)Y_\ell^m, $$
so that by orthonormality
$$ \sum_{\ell=0}^L\sum_{m=-\ell}^{\ell} |\widehat f(\ell,m)|^2 = \|f\|_{L^2(S^2)}^2. $$
The constant spherical harmonic is $Y_0^0=\frac{1}{\sqrt{4\pi}}$, and hence
$$ |\widehat f(0,0)| = \frac{1}{\sqrt{4\pi}} \left|\int_{S^2}f(\omega)d\sigma(\omega)\right|. $$
For $\ell\geq1$ we insert the eigenvalue $\lambda_\ell=\ell(\ell+1)$ and apply Cauchy--Schwarz to obtain
$$ \begin{aligned} \sum_{\ell=1}^L\sum_{m=-\ell}^{\ell}|\widehat f(\ell,m)| &= \sum_{\ell=1}^L\sum_{m=-\ell}^{\ell} \frac{1}{\sqrt{\lambda_\ell}} \sqrt{\lambda_\ell}|\widehat f(\ell,m)|\\ &\leq \left( \sum_{\ell=1}^L\sum_{m=-\ell}^{\ell}\frac{1}{\lambda_\ell} \right)^{\frac{1}{2}} \left( \sum_{\ell=1}^L\sum_{m=-\ell}^{\ell} \lambda_\ell|\widehat f(\ell,m)|^2 \right)^{\frac{1}{2}}. \end{aligned} $$
The sum in the second factor is the spherical Dirichlet energy, which is to say that
$$ \sum_{\ell=1}^L\sum_{m=-\ell}^{\ell} \lambda_\ell|\widehat f(\ell,m)|^2 = \|\nabla_{S^2}f\|_{L^2(S^2)}^2. $$
As for the first factor, the multiplicity of the eigenvalue $\lambda_\ell$ is $2\ell+1$, so
$$ \sum_{\ell=1}^L\sum_{m=-\ell}^{\ell}\frac{1}{\ell(\ell+1)} = \sum_{\ell=1}^L\frac{2\ell+1}{\ell(\ell+1)} = \sum_{\ell=1}^L\left(\frac{1}{\ell}+\frac{1}{\ell+1}\right) \leq 2(1+\log L). $$
Combining these identities and dividing by $\|f\|_{L^2(S^2)}$ we obtain the first estimate. Since $S^2$ has surface area $4\pi$, we have $\|\nabla_{S^2}f\|_{L^2(S^2)}\leq\sqrt{4\pi}\|f\|_{C^1(S^2)}$, which implies the second estimate. Finally, if $f\geq 0$, then $\mu\leq\sqrt{4\pi}\|f\|_{L^2(S^2)}$ by Cauchy--Schwarz, and using this inequality in both terms we arrive at the final assertion.
\end{proof}

\section{Sampling and reconstruction}\label{sec:recovery}

We now combine the analytic bounds of the previous two sections with the effect of moving the sensors and with the recovery result from Section~\ref{sec:tools}. A displacement of a sensor is first converted into a measurement error, and the Fourier ratio then bounds the approximation term in stable $\ell^1$ recovery. Partitions into cells of equal volume provide a separate description of the number of random sensors needed to occupy a prescribed number of cells. Whenever a ratio of norms of $f$ occurs in this section, we assume that $f$ is nonzero.

\subsection{Sensor placement and random sensors on the square}

We begin by allowing the sampling location in each cell to vary independently. Such perturbations do not preserve the Fourier structure exactly, so it is better to regard them as measurement errors relative to the uniform grid.

\begin{theorem}\label{thm:square-perturbation}
Let $f\in C^1([0,1]^2)$ be $1$-periodic, and define $g(x)=f(\frac{x}{N})$ for $x\in\mathbb Z_N^2$. Let $0\leq\delta\leq1$ and let $1\leq q\leq N^2$. For each $x\in\mathbb Z_N^2$, choose $t_x\in[0,\delta]^2$ and define $g_t(x)=f(\frac{x+t_x}{N})$. Then
$$ |g_t(x)-g(x)| \leq \frac{2\delta}{N}\|f\|_{C^1([0,1]^2)}. $$
Consequently, for every set $X\subset\mathbb Z_N^2$ of cardinality $q$,
$$ \left( \frac{1}{q}\sum_{x\in X}|g_t(x)-g(x)|^2 \right)^{\frac{1}{2}} \leq \frac{2\delta}{N}\|f\|_{C^1([0,1]^2)}. $$
In particular, this error is at most $\varepsilon\|f\|_{L^2([0,1]^2)}$ whenever
$$ \delta \leq \frac{\varepsilon N}{2} \frac{\|f\|_{L^2([0,1]^2)}}{\|f\|_{C^1([0,1]^2)}}. $$
\end{theorem}

\begin{proof}
By the fundamental theorem of calculus,
$$ g_t(x)-g(x) = \int_0^1 \nabla f\left(\frac{x+st_x}{N}\right) \cdot\frac{t_x}{N}ds, $$
and since both coordinates of $t_x$ are at most $\delta$,
$$ |g_t(x)-g(x)| \leq \frac{\delta}{N} \left( \left\|\frac{\partial f}{\partial x_1}\right\|_{L^\infty} + \left\|\frac{\partial f}{\partial x_2}\right\|_{L^\infty} \right), $$
which proves the pointwise estimate. The averaged estimate and the final assertion follow at once.
\end{proof}

\begin{corollary}\label{cor:full-cell}
If
$$ N \geq \frac{2}{\varepsilon} \frac{\|f\|_{C^1([0,1]^2)}}{\|f\|_{L^2([0,1]^2)}}, $$
then the conclusion of Theorem~\ref{thm:square-perturbation} holds with error at most $\varepsilon\|f\|_{L^2([0,1]^2)}$ for every choice $t_x\in[0,1)^2$.
\end{corollary}

We complete our discussion of placement on the square by describing what happens when the sensors themselves are placed uniformly at random.

\begin{theorem}\label{thm:coupon}
Let $f\in C^1([0,1]^2)$ be periodic, and let $Y_1,Y_2,\ldots$ be independent uniformly distributed points in $[0,1)^2$. Let $Z_j\in\mathbb Z_N^2$ be the unique index such that $Y_j\in Q_{Z_j}$, where $Q_x=\frac{x}{N}+[0,N^{-1})^2$ is the grid cell indexed by $x$, and set $X_m=\{Z_1,\ldots,Z_m\}$. For $1\leq q\leq N^2$, define
$$ T_q=\min\{m\geq 1:|X_m|=q\}. $$
Then
$$ \mathbb E[T_q] = N^2(H_{N^2}-H_{N^2-q}), $$
where $H_0=0$ and $H_k=\sum_{j=1}^k\frac{1}{j}$. Moreover, $X_{T_q}$ is uniformly distributed among the $q$-element subsets of $\mathbb Z_N^2$.

For each $x\in X_{T_q}$, use the first sensor that lands in $Q_x$, and denote its location by $Y_x$. If $g_t(x)=f(Y_x)$ and $g(x)=f(\frac{x}{N})$, then
$$ \left( \frac{1}{q}\sum_{x\in X_{T_q}}|g_t(x)-g(x)|^2 \right)^{\frac{1}{2}} \leq \frac{2}{N}\|f\|_{C^1([0,1]^2)}. $$
\end{theorem}

\begin{proof}
Every cell has area $\frac{1}{N^2}$, and since the points $Y_j$ are independent and uniformly distributed, the cell indices $Z_j$ are independent and uniform on $\mathbb Z_N^2$. Suppose that $k$ distinct cells have already been hit. There are $N^2-k$ unoccupied cells, so the probability that the next sensor lands in a new cell is $p_k=\frac{N^2-k}{N^2}$. Let $W_k$ be the number of additional sensors needed to move from $k$ occupied cells to $k+1$. Conditional on the cells already encountered, $W_k$ has the geometric distribution with success probability $p_k$, and therefore
$$ \mathbb E[W_k]=\frac{1}{p_k}=\frac{N^2}{N^2-k}. $$
Since $T_q=\sum_{k=0}^{q-1}W_k$, linearity of expectation implies that
$$ \mathbb E[T_q] = \sum_{k=0}^{q-1}\frac{N^2}{N^2-k} = N^2\sum_{j=N^2-q+1}^{N^2}\frac{1}{j} = N^2(H_{N^2}-H_{N^2-q}). $$

We now turn to the uniformity assertion. Let $\pi$ be any permutation of the $N^2$ cells. The sequence $\pi(Z_1),\pi(Z_2),\ldots$ has the same distribution as the original sequence, and if $A$ and $B$ are two $q$-element subsets, one can choose $\pi$ with $\pi(A)=B$. It follows that $\mathbb P(X_{T_q}=A)=\mathbb P(X_{T_q}=B)$, and since the probabilities over all $\binom{N^2}{q}$ possible sets sum to one,
$$ \mathbb P(X_{T_q}=A)=\frac{1}{\binom{N^2}{q}}. $$
Finally, since $Y_x\in Q_x$, there is a vector $t_x\in[0,1)^2$ such that $Y_x=\frac{x+t_x}{N}$, namely the one with coordinates $t_{x,j}=NY_{x,j}-x_j$, and the last estimate is therefore Theorem~\ref{thm:square-perturbation} with $\delta=1$.
\end{proof}

\begin{remark}
Theorem~\ref{thm:coupon} provides the uniform random subset required by many discrete recovery theorems, together with a deterministic bound for the perturbation of the measurements. A particular recovery theorem may be applied after checking its requirement on the sample size and matching its data fidelity norm to the error bound above. We emphasize that the coupon collector argument does not by itself assert any recovery probability.
\end{remark}

\subsection{Recovery on the periodic square}

We now put together the Fourier ratio estimate on the square, the placement bound, and Proposition~\ref{prop:bos-recovery}. The theorem below uses independent grid indices sampled with replacement, whereas Theorem~\ref{thm:coupon} describes a different model, in which one retains the first sensor in each of a prescribed number of distinct cells.

\begin{theorem}\label{thm:square-random-recovery}
Let $f\in C^1([0,1]^2)$ be nonzero and $1$-periodic, let $N\geq2$, and fix $t\in[0,1)^2$. Define $g_t(x)=f(\frac{x+t}{N})$ for $x\in\mathbb Z_N^2$, assume that this sampled vector is nonzero, and suppose that $FR(g_t)\leq r$. Fix $0<\varepsilon<\frac{1}{2}$ and $\gamma\geq1$, and set
$$ s= \min\left\{ N^2, \left\lceil\frac{r^2}{\varepsilon^2}\right\rceil \right\}. $$
Choose $X_1,\ldots,X_q$ independently and uniformly from $\mathbb Z_N^2$. For each $i$, choose a vector $u_i\in[0,1)^2$ and observe
$$ y_i = f\left(\frac{X_i+t+u_i}{N}\right), $$
where the periodic extension of $f$ is used if a point crosses the boundary of the unit square. Let $g^*:\mathbb Z_N^2\to\mathbb C$ minimize $\|\widehat h\|_1$ among all functions $h:\mathbb Z_N^2\to\mathbb C$ satisfying
$$ \left( \frac{1}{q}\sum_{i=1}^q|h(X_i)-y_i|^2 \right)^{\frac{1}{2}} \leq \frac{2}{N}\|f\|_{C^1([0,1]^2)}. $$
There is a constant $C_\gamma$ such that, if
$$ q \geq C_\gamma s\log^2(2s)\log(2N^2), $$
then, with probability at least $1-\frac{1}{N^{2\gamma}}$,
$$ \frac{1}{N}\|g^*-g_t\|_2 \leq \frac{3\varepsilon}{N}\|g_t\|_2 + \frac{10}{N}\|f\|_{C^1([0,1]^2)}. $$
The probability is taken over the choice of the grid indices $X_i$, and the conclusion holds simultaneously for every choice of the displacement vectors $u_i$, including choices that depend on the sampled indices.
\end{theorem}

\begin{proof}
Equip $\mathbb Z_N^2$ with normalized counting measure and consider the characters $\phi_m(x)=e^{\frac{2\pi i x\cdot m}{N}}$, $m\in\mathbb Z_N^2$. These characters form an orthonormal basis and satisfy $|\phi_m(x)|=1$, so Proposition~\ref{prop:bos-recovery} applies with $D=N^2$ and $K=1$. The coefficient of $g_t$ relative to normalized counting measure is
$$ a(m) = \frac{1}{N^2}\sum_{x\in\mathbb Z_N^2} e^{-\frac{2\pi i x\cdot m}{N}}g_t(x) = \frac{1}{N}\widehat{g_t}(m), $$
and consequently
$$ \frac{\|a\|_1}{\|a\|_2}=FR(g_t)\leq r \text{ and } \|a\|_2=\frac{1}{N}\|g_t\|_2. $$
The displacement of the $i$th sensor from its reference point is $\frac{u_i}{N}$, and the argument with the fundamental theorem of calculus used in the proof of Theorem~\ref{thm:square-perturbation} shows that
$$ |y_i-g_t(X_i)| \leq \frac{2}{N}\|f\|_{C^1([0,1]^2)}. $$
Define the sampling matrix $A_{i,m}=\frac{1}{\sqrt q}\phi_m(X_i)$ and normalize the observations by setting $\widetilde y_i=\frac{1}{\sqrt q}y_i$. Then $\widetilde y=Aa+e$, where
$$ \|e\|_2 \leq \frac{2}{N}\|f\|_{C^1([0,1]^2)}. $$
If $s<N^2$, then the definition of $s$ implies that
$$ \frac{\sigma_s(a)_1}{\sqrt s} \leq \frac{\|a\|_1}{\sqrt s} \leq \varepsilon\|a\|_2, $$
and if $s=N^2$, then $\sigma_s(a)_1=0$, so the same conclusion holds. By Proposition~\ref{prop:bos-recovery},
$$ \|a^*-a\|_2 \leq 3\varepsilon\|a\|_2 + \frac{10}{N}\|f\|_{C^1([0,1]^2)}. $$
The coefficient vector of a candidate $h$ is $\frac{1}{N}\widehat h$, and multiplying the objective by $\frac{1}{N}$ does not change its minimizer, so $a^*$ is the coefficient vector of $g^*$. Since $\|a^*-a\|_2=\frac{1}{N}\|g^*-g_t\|_2$ by Parseval's identity for normalized counting measure, the proof is complete.
\end{proof}

\begin{corollary}\label{cor:translated-random-recovery}
Let $f\in C^1([0,1]^2)$ be nonzero and $1$-periodic. There is a translation $t\in[0,1)^2$ for which Theorem~\ref{thm:square-random-recovery} applies with
$$ r = 1+ \frac{1}{\sqrt{2}}\sqrt{1+\log N} \frac{\|\nabla f\|_{L^2([0,1]^2)}}{\|f\|_{L^2([0,1]^2)}}. $$
If the $u_i$ are chosen independently and uniformly from $[0,1)^2$, independently of the grid indices, then the physical sensor locations $\frac{X_i+t+u_i}{N}$ are independent and uniformly distributed on the two-dimensional torus.
\end{corollary}

\begin{proof}
Choose the translation provided by Theorem~\ref{thm:translated-grid}. Its Fourier ratio is at most the value of $r$ in the statement, so Theorem~\ref{thm:square-random-recovery} applies. For the last assertion, the sensor locations are interpreted modulo $\mathbb Z^2$. The cells $\frac{x+t}{N}+[0,N^{-1})^2$, indexed by all $x\in\mathbb Z_N^2$ and interpreted periodically, partition the torus and have equal area, so choosing a cell uniformly and then choosing a point uniformly inside it produces normalized Lebesgue measure on the whole torus.
\end{proof}

\begin{corollary}[Piecewise constant reconstruction]\label{cor:square-continuous-recovery}
Under the hypotheses of Theorem~\ref{thm:square-random-recovery}, let
$$ Q_{x,t}=\frac{x+t}{N}+[0,N^{-1})^2 \text{ for } x\in\mathbb Z_N^2, $$
with addition understood periodically, and define $F^*$ to equal $g^*(x)$ on $Q_{x,t}$. On the recovery event in that theorem,
$$ \|F^*-f\|_{L^2([0,1]^2)} \leq3\varepsilon\|f\|_{L^2([0,1]^2)} +\frac{12+6\varepsilon}{N}\|f\|_{C^1([0,1]^2)}. $$
\end{corollary}

\begin{proof}
Let $G_t$ equal $g_t(x)$ on $Q_{x,t}$. Each cell has area $N^{-2}$, so
$$ \|F^*-G_t\|_{L^2}=\frac{1}{N}\|g^*-g_t\|_2 \text{ and } \|G_t\|_{L^2}=\frac{1}{N}\|g_t\|_2, $$
and the fundamental theorem of calculus inside each cell shows that $\|G_t-f\|_{L^2}\leq\frac{2}{N}\|f\|_{C^1}$. We now apply the triangle inequality and Theorem~\ref{thm:square-random-recovery}. The approximation term satisfies
$$ \frac{3\varepsilon}{N}\|g_t\|_2 \leq3\varepsilon\|f\|_{L^2} +\frac{6\varepsilon}{N}\|f\|_{C^1}, $$
and the result follows.
\end{proof}

\subsection{Sensor stability on compact manifolds}

We now turn to compact manifolds. Write $d_g$ for the geodesic distance on $M$, and for $f\in C^1(M)$ let
$$ \|f\|_{C^1(M)} =\max\{\|f\|_{L^\infty(M)}, \|\nabla_gf\|_{L^\infty(M)}\}. $$
The statements below apply to real or complex valued functions.

\begin{proposition}\label{prop:manifold-placement-source}
Let $(M,g)$ be a connected compact Riemannian manifold without boundary, and let $f$ be a nonzero function in $C^1(M)$. Let $q\geq1$, let $x_1,\ldots,x_q$ be reference sensor locations and $y_1,\ldots,y_q$ perturbed sensor locations satisfying $d_g(x_i,y_i)\leq\delta$ for $i=1,\ldots,q$, and define $g(x_i)=f(x_i)$ and $g_t(x_i)=f(y_i)$. Then
$$ |g_t(x_i)-g(x_i)|\leq\delta\|f\|_{C^1(M)} $$
for every $i$, and
$$ \left( \frac{|M|}{q}\sum_{i=1}^q|g_t(x_i)-g(x_i)|^2 \right)^{\frac{1}{2}} \leq |M|^{\frac{1}{2}}\delta\|f\|_{C^1(M)}. $$
\end{proposition}

\begin{proof}
Join $x_i$ to $y_i$ by a minimizing geodesic $\gamma_i:[0,\ell_i]\to M$ parameterized by arclength. Since $\ell_i=d_g(x_i,y_i)\leq\delta$, the fundamental theorem of calculus along $\gamma_i$ shows that
$$ f(y_i)-f(x_i) = \int_0^{\ell_i} \langle\nabla_gf(\gamma_i(s)),\gamma_i'(s)\rangle_gds, $$
and since $|\gamma_i'(s)|_g=1$,
$$ |f(y_i)-f(x_i)| \leq \ell_i\|\nabla_gf\|_{L^\infty(M)} \leq \delta\|f\|_{C^1(M)}. $$
The second assertion follows by squaring this estimate, summing over $i$, multiplying by $\frac{|M|}{q}$, and taking the square root.
\end{proof}

This estimate can be read as a tolerance condition on how far the actual sensors may move if the total measurement error is to be kept at a prescribed relative scale.

\begin{corollary}\label{cor:manifold-placement-tolerance}
Under the hypotheses of Proposition~\ref{prop:manifold-placement-source}, if
$$ \delta \leq \frac{\varepsilon\|f\|_{L^2(M)}}{|M|^{\frac{1}{2}}\|f\|_{C^1(M)}}, $$
then
$$ \left( \frac{|M|}{q}\sum_{i=1}^q|g_t(x_i)-g(x_i)|^2 \right)^{\frac{1}{2}} \leq \varepsilon\|f\|_{L^2(M)}. $$
\end{corollary}

\begin{proof}
This follows by substituting the assumed upper bound for $\delta$ into Proposition~\ref{prop:manifold-placement-source}.
\end{proof}

If a sensor is only known to lie in the same cell as its intended location, then the diameter of the cell can be used as the displacement parameter.

\begin{corollary}\label{cor:manifold-partition-cells}
Let $(M,g)$, $f$, $x_i$, $y_i$, $g$, and $g_t$ be as in Proposition~\ref{prop:manifold-placement-source}. Suppose that $x_i,y_i\in Q_i\subset M$ and $\operatorname{diam}_g(Q_i)\leq h$ for every $i$. Then $d_g(x_i,y_i)\leq h$. Consequently, if
$$ h \leq \frac{\varepsilon\|f\|_{L^2(M)}}{|M|^{\frac{1}{2}}\|f\|_{C^1(M)}}, $$
then
$$ \left( \frac{|M|}{q}\sum_{i=1}^q|g_t(x_i)-g(x_i)|^2 \right)^{\frac{1}{2}} \leq \varepsilon\|f\|_{L^2(M)}. $$
\end{corollary}

\begin{proof}
Since $x_i,y_i\in Q_i$, we have $d_g(x_i,y_i)\leq\operatorname{diam}_g(Q_i)\leq h$, and the conclusion follows from Proposition~\ref{prop:manifold-placement-source} with $\delta=h$ and the assumed upper bound for $h$.
\end{proof}

If the cells come from a partition into $K$ pieces whose diameters decay at the natural rate $K^{-\frac{1}{d}}$, then the same condition can be stated as a lower bound on the number of cells.

\begin{corollary}\label{cor:manifold-cell-count}
Assume the hypotheses of Corollary~\ref{cor:manifold-partition-cells}. Suppose further that $M$ has dimension $d$ and that there exist $A>0$ and a positive integer $K$ such that $\operatorname{diam}_g(Q_i)\leq AK^{-\frac{1}{d}}$ for every $i$. If
$$ K \geq \left( \frac{A|M|^{\frac{1}{2}}\|f\|_{C^1(M)}}{\varepsilon\|f\|_{L^2(M)}} \right)^d, $$
then
$$ \left( \frac{|M|}{q}\sum_{i=1}^q|g_t(x_i)-g(x_i)|^2 \right)^{\frac{1}{2}} \leq \varepsilon\|f\|_{L^2(M)}. $$
\end{corollary}

\begin{proof}
The assumed lower bound for $K$ implies that
$$ AK^{-\frac{1}{d}} \leq \frac{\varepsilon\|f\|_{L^2(M)}}{|M|^{\frac{1}{2}}\|f\|_{C^1(M)}}, $$
and the conclusion follows from Corollary~\ref{cor:manifold-partition-cells}.
\end{proof}

Just as on the square, partitions into cells of equal volume allow us to pass from uniformly random physical sensors to uniformly random subsets of reference cells.

\begin{theorem}\label{thm:manifold-coupon-source}
Let $(M,g)$ be a connected compact Riemannian manifold without boundary, and let $Q_1,\ldots,Q_K$ be a measurable partition of $M$ satisfying $|Q_k|=\frac{|M|}{K}$ for every $k$. Choose a representative point $x_k\in Q_k$. Let $Y_1,Y_2,\ldots$ be independent uniformly distributed random sensor locations on $M$ with respect to normalized Riemannian volume, and let $Z_j$ be the unique index such that $Y_j\in Q_{Z_j}$. For $1\leq q\leq K$, define
$$ X_m=\{Z_1,\ldots,Z_m\} \text{ and } T_q=\min\{m\geq1:|X_m|=q\}. $$
Then
$$ \mathbb E[T_q]=K(H_K-H_{K-q}), $$
where $H_0=0$ and $H_k=\sum_{j=1}^k\frac{1}{j}$, and $X_{T_q}$ is uniformly distributed among the $q$-element subsets of $\{1,\ldots,K\}$.

Suppose further that $\operatorname{diam}_g(Q_k)\leq h_K$ for every $k$. Fix $1\leq q_{\mathrm{rec}}\leq K$, write $T_{\mathrm{rec}}=T_{q_{\mathrm{rec}}}$, and, for each $k\in X_{T_{\mathrm{rec}}}$, let $Y_k^*$ denote the first sensor to land in $Q_k$. If $f\in C^1(M)$, $g(k)=f(x_k)$, and $g_t(k)=f(Y_k^*)$, then
$$ \left( \frac{|M|}{q_{\mathrm{rec}}} \sum_{k\in X_{T_{\mathrm{rec}}}} |g_t(k)-g(k)|^2 \right)^{\frac{1}{2}} \leq |M|^{\frac{1}{2}}h_K\|f\|_{C^1(M)}. $$
In particular, the placement error is at most $\varepsilon\|f\|_{L^2(M)}$ whenever
$$ h_K \leq \frac{\varepsilon\|f\|_{L^2(M)}}{|M|^{\frac{1}{2}}\|f\|_{C^1(M)}}. $$
\end{theorem}

\begin{proof}
The hypothesis that the cells have equal volume implies that
$$ \mathbb P(Z_j=k) = \frac{|Q_k|}{|M|} = \frac{1}{K}, $$
so the cell indices are independent and uniformly distributed. The coupon collector calculation in the proof of Theorem~\ref{thm:coupon} then shows that $\mathbb E[T_q] = K(H_K-H_{K-q})$, and invariance under permutations of the cell labels shows that $X_{T_q}$ is uniformly distributed among the $q$-element subsets. For $k\in X_{T_{\mathrm{rec}}}$, both $x_k$ and $Y_k^*$ belong to $Q_k$, so $d_g(x_k,Y_k^*)\leq\operatorname{diam}_g(Q_k)\leq h_K$, and Proposition~\ref{prop:manifold-placement-source} implies that
$$ |g_t(k)-g(k)| \leq h_K\|f\|_{C^1(M)}. $$
The placement estimate follows by squaring, summing over the occupied cells, multiplying by $\frac{|M|}{q_{\mathrm{rec}}}$, and taking the square root, and the final assertion follows from the condition on $h_K$.
\end{proof}

\subsection{Sensors on the sphere}

On the sphere the geodesic placement estimate takes a particularly simple form, and a partition into cells of equal area leads to the same coupon collector description of random physical sensors as on the square. We write $d_{S^2}$ for the geodesic distance on $S^2$.

\begin{proposition}\label{prop:sphere-perturbation}
Let $f\in C^1(S^2)$ and $q\geq1$, and let $\omega_1,\ldots,\omega_q$ and $\omega_{1,t},\ldots,\omega_{q,t}$ be two collections of sensor locations satisfying $d_{S^2}(\omega_{i,t},\omega_i)\leq\delta$ for every $i$. Then
$$ |f(\omega_{i,t})-f(\omega_i)| \leq \delta\|f\|_{C^1(S^2)} $$
and
$$ \left( \frac{4\pi}{q}\sum_{i=1}^q |f(\omega_{i,t})-f(\omega_i)|^2 \right)^{\frac{1}{2}} \leq \sqrt{4\pi}\delta\|f\|_{C^1(S^2)}. $$
Thus the relative placement error is bounded by
$$ \sqrt{4\pi}\delta \frac{\|f\|_{C^1(S^2)}}{\|f\|_{L^2(S^2)}}. $$
\end{proposition}

\begin{proof}
This is Proposition~\ref{prop:manifold-placement-source} with $M=S^2$ and $|M|=4\pi$.
\end{proof}

\begin{corollary}\label{cor:sphere-cells}
Suppose $Q_1,\ldots,Q_q\subset S^2$ satisfy $\operatorname{diam}_{S^2}(Q_i)\leq\frac{A}{N}$. If $\omega_i,\omega_{i,t}\in Q_i$, then $d_{S^2}(\omega_{i,t},\omega_i)\leq\frac{A}{N}$. Consequently, the averaged placement error in Proposition~\ref{prop:sphere-perturbation} is at most $\varepsilon\|f\|_{L^2(S^2)}$ whenever
$$ N \geq \frac{A\sqrt{4\pi}}{\varepsilon} \frac{\|f\|_{C^1(S^2)}}{\|f\|_{L^2(S^2)}}. $$
\end{corollary}

\begin{theorem}\label{thm:sphere-coupon}
Let $f\in C^1(S^2)$, and let $Q_1,\ldots,Q_K$ be a measurable partition of $S^2$ such that $\sigma(Q_k)=\frac{4\pi}{K}$ for every $k=1,\ldots,K$, and choose a representative point $\omega_k\in Q_k$ for each cell. Let $Y_1,Y_2,\ldots$ be independent uniformly distributed random sensor locations on $S^2$. For $j\geq1$, let $Z_j\in\{1,\ldots,K\}$ be the unique cell index such that $Y_j\in Q_{Z_j}$, and define $X_m=\{Z_1,\ldots,Z_m\}$. For $1\leq q\leq K$, set $T_q=\min\{m\geq1:|X_m|=q\}$. Then
$$ \mathbb E[T_q]=K(H_K-H_{K-q}), $$
where $H_0=0$ and $H_k=\sum_{j=1}^k\frac{1}{j}$. Moreover, $X_{T_q}$ is uniformly distributed among the $q$-element subsets of $\{1,\ldots,K\}$.

Suppose further that $\operatorname{diam}_{S^2}(Q_k)\leq h_K$ for every $k$. Fix $1\leq q_{\mathrm{rec}}\leq K$, write $T_{\mathrm{rec}}=T_{q_{\mathrm{rec}}}$, and, for each $k\in X_{T_{\mathrm{rec}}}$, let $Y_k^*$ denote the first sensor location to land in $Q_k$. If $g(k)=f(\omega_k)$ and $g_t(k)=f(Y_k^*)$, then
$$ \left( \frac{4\pi}{q_{\mathrm{rec}}} \sum_{k\in X_{T_{\mathrm{rec}}}} |g_t(k)-g(k)|^2 \right)^{\frac{1}{2}} \leq \sqrt{4\pi}h_K\|f\|_{C^1(S^2)}. $$
In particular, this placement error is at most $\varepsilon\|f\|_{L^2(S^2)}$ whenever
$$ h_K \leq \frac{\varepsilon\|f\|_{L^2(S^2)}}{\sqrt{4\pi}\|f\|_{C^1(S^2)}}. $$
\end{theorem}

\begin{proof}
This follows by applying Theorem~\ref{thm:manifold-coupon-source} to the given partition of $S^2$, since the normalized surface measure assigns probability $\frac{1}{K}$ to each cell and $|S^2|=4\pi$.
\end{proof}

\subsection{Recovery in the Laplace--Beltrami eigenbasis}

The results above control the sensor locations themselves. We now return to the Laplace--Beltrami eigenbasis and combine this control with the Fourier ratio bounds of Section~\ref{sec:manifolds} and with Proposition~\ref{prop:bos-recovery}. For the spectral space $V_\Lambda$, let
$$ K_\Lambda = |M|^{\frac{1}{2}} \sup_{\lambda_j\leq \Lambda}\|e_j\|_{L^\infty(M)}. $$
The functions $\phi_j=|M|^{\frac{1}{2}}e_j$ form an orthonormal system in $L^2(M,\frac{1}{|M|}dV_g)$, and $K_\Lambda$ is the boundedness parameter of this normalized system.

\begin{theorem}\label{thm:manifold-recovery}
Let $f\in V_\Lambda$ be nonzero and suppose that $FR_\Lambda(f)\leq r$. Fix $0<\varepsilon<\frac{1}{2}$ and $\gamma\geq 1$, and set
$$ s= \min\left\{ D_\Lambda, \left\lceil\frac{r^2}{\varepsilon^2}\right\rceil \right\}. $$
Let $x_1,\ldots,x_q$ be independent random points chosen according to the normalized Riemannian measure. Suppose that the observations have the form $y_i=f(x_i)+\xi_i$ and satisfy
$$ \left( \frac{|M|}{q}\sum_{i=1}^q|\xi_i|^2 \right)^{\frac{1}{2}} \leq \eta. $$
Let $f^*\in V_\Lambda$ minimize $\|\widehat h\|_1$ among all $h\in V_\Lambda$ satisfying
$$ \left( \frac{|M|}{q}\sum_{i=1}^q|h(x_i)-y_i|^2 \right)^{\frac{1}{2}} \leq \eta. $$
There is a constant $C_\gamma$ such that, if
$$ q \geq C_\gamma K_\Lambda^2s\log^2(2sK_\Lambda^2)\log(2D_\Lambda), $$
then, with probability at least $1-\frac{1}{D_\Lambda^\gamma}$,
$$ \|f^*-f\|_{L^2(M)} \leq 3\varepsilon\|f\|_{L^2(M)}+5\eta. $$
\end{theorem}

\begin{proof}
Write $f=\sum_{\lambda_j\leq \Lambda}a_j\phi_j$. Since $\phi_j=|M|^{\frac{1}{2}}e_j$, comparison with the original eigenfunction expansion shows that $a_j=\frac{1}{\sqrt{|M|}}\widehat f(j)$. The common factor cancels from the ratio of the coefficient norms, so the coefficient vector satisfies
$$ \frac{\|a\|_1}{\|a\|_2} =FR_\Lambda(f)\leq r. $$
If $s<D_\Lambda$, the definition of $s$ and the elementary bound $\sigma_s(a)_1\leq\|a\|_1$ imply that
$$ \frac{\sigma_s(a)_1}{\sqrt{s}} \leq \frac{\|a\|_1}{\sqrt{s}} \leq \varepsilon\|a\|_2, $$
and if $s=D_\Lambda$, then $\sigma_s(a)_1=0$, so the same conclusion holds.

Consider the sampling matrix $A_{ij}=\frac{1}{\sqrt q}\phi_j(x_i)$, so that the vector $Aa$ has entries $\frac{1}{\sqrt q}f(x_i)$. We normalize the observations and the noise by setting $\widetilde y_i=\frac{1}{\sqrt q}y_i$ and $e_i=\frac{1}{\sqrt q}\xi_i$. Then $\widetilde y=Aa+e$, and the hypothesis on the noise in physical space implies that
$$ \|e\|_2 = \left(\frac{1}{q}\sum_{i=1}^q|\xi_i|^2\right)^{\frac{1}{2}} \leq \frac{1}{\sqrt{|M|}}\eta. $$
The data fidelity constraint in the theorem is therefore exactly the constraint in Proposition~\ref{prop:bos-recovery}, with noise tolerance $\frac{1}{\sqrt{|M|}}\eta$. Moreover, $\|\widehat h\|_1=|M|^{\frac{1}{2}}\|b\|_1$, where $b$ is the coefficient vector of $h$ in the $\phi_j$ basis, and multiplication of the objective by the fixed positive number $|M|^{\frac{1}{2}}$ does not change its minimizer. It follows from Proposition~\ref{prop:bos-recovery} that
$$ \|a^*-a\|_2 \leq 3\frac{\sigma_s(a)_1}{\sqrt{s}} + \frac{5}{\sqrt{|M|}}\eta, $$
where $a^*$ is the coefficient vector of $f^*$. Since the $\phi_j$ are orthonormal for the normalized measure,
$$ \|a\|_2 = \frac{1}{\sqrt{|M|}}\|f\|_{L^2(M)}, $$
and the same scaling relates $\|a^*-a\|_2$ to $\|f^*-f\|_{L^2(M)}$. The theorem follows by substituting the approximation estimate.
\end{proof}

\begin{corollary}[Explicit relative error]\label{cor:manifold-recovery-explicit}
Let $f\in V_\Lambda$ be nonzero and suppose $FR_\Lambda(f)\leq r$. Fix $0<\varepsilon<\frac{1}{2}$ and $\gamma\geq1$, and set
$$ s=\min\left\{D_\Lambda,\left\lceil\frac{r^2}{\varepsilon^2}\right\rceil\right\}. $$
Choose $x_1,\ldots,x_q$ independently according to normalized Riemannian measure, and let $f^*\in V_\Lambda$ minimize $\|\widehat h\|_1$ subject to
$$ \left(\frac{|M|}{q}\sum_{i=1}^q|h(x_i)-f(x_i)|^2\right)^{\frac{1}{2}} \leq\varepsilon\|f\|_{L^2(M)}. $$
There is a constant $C_\gamma$ such that, if
$$ q\geq C_\gamma K_\Lambda^2s\log^2(2sK_\Lambda^2)\log(2D_\Lambda), $$
then, with probability at least $1-\frac{1}{D_\Lambda^\gamma}$,
$$ \|f^*-f\|_{L^2(M)}\leq8\varepsilon\|f\|_{L^2(M)}. $$
\end{corollary}

\begin{proof}
We apply Theorem~\ref{thm:manifold-recovery} with $\xi_i=0$ and $\eta=\varepsilon\|f\|_{L^2(M)}$. The zero noise satisfies this tolerance, the reconstruction constraint is the one in the statement, and the two error terms in the theorem are bounded by $3\varepsilon\|f\|_{L^2(M)}$ and $5\varepsilon\|f\|_{L^2(M)}$, respectively.
\end{proof}

\begin{corollary}\label{cor:manifold-perturbed-recovery}
Assume the hypotheses of Theorem~\ref{thm:manifold-recovery}. After choosing the reference points $x_1,\ldots,x_q$, suppose the actual sensor locations $x_{1,t},\ldots,x_{q,t}$ satisfy $d_g(x_{i,t},x_i)\leq\delta$ and that the observations are $y_i=f(x_{i,t})$. Use the data fidelity tolerance
$$ \eta = \sqrt{|M|}\delta\|\nabla_gf\|_{L^\infty(M)}. $$
If
$$ q \geq C_\gamma K_\Lambda^2s\log^2(2sK_\Lambda^2)\log(2D_\Lambda), $$
then, with probability at least $1-\frac{1}{D_\Lambda^\gamma}$,
$$ \|f^*-f\|_{L^2(M)} \leq 3\varepsilon\|f\|_{L^2(M)} + 5\sqrt{|M|}\delta\|\nabla_gf\|_{L^\infty(M)}. $$
The conclusion holds simultaneously for every choice of the displaced points satisfying the stated distance bound.
\end{corollary}

\begin{proof}
Joining $x_i$ to $x_{i,t}$ by a minimizing geodesic and applying the fundamental theorem of calculus along that geodesic, we see that
$$ |f(x_{i,t})-f(x_i)| \leq \delta\|\nabla_gf\|_{L^\infty(M)}, $$
and consequently
$$ \left( \frac{|M|}{q}\sum_{i=1}^q|f(x_{i,t})-f(x_i)|^2 \right)^{\frac{1}{2}} \leq \sqrt{|M|}\delta\|\nabla_gf\|_{L^\infty(M)}. $$
The conclusion is therefore exactly Theorem~\ref{thm:manifold-recovery} with this value of $\eta$.
\end{proof}

\begin{corollary}\label{cor:manifold-recovery}
Let $f\in V_\Lambda$ be nonzero. For $s_0>0$, define
$$ r_{\Lambda,s_0}(f) = \frac{1}{\sqrt{|M|}} \frac{\left|\int_Mf(x)dV_g(x)\right|}{\|f\|_{L^2(M)}} + \left(\sum_{0<\lambda_j\leq \Lambda}\frac{1}{\lambda_j^{s_0}}\right)^{\frac{1}{2}} \frac{\|(-\Delta_g)^{\frac{s_0}{2}}f\|_{L^2(M)}}{\|f\|_{L^2(M)}}. $$
Then $FR_\Lambda(f)\leq r_{\Lambda,s_0}(f)$. Consequently, Theorem~\ref{thm:manifold-recovery} and Corollary~\ref{cor:manifold-recovery-explicit} apply with $r=r_{\Lambda,s_0}(f)$.
\end{corollary}

\begin{proof}
The inequality $FR_\Lambda(f)\leq r_{\Lambda,s_0}(f)$ is the first estimate of Proposition~\ref{prop:manifold-fr} with $s=s_0$.
\end{proof}

\begin{remark}
The factor $K_\Lambda^2$ measures the concentration of the selected eigenfunctions through their supremum norms, and on manifolds where eigenfunctions concentrate strongly the resulting sample bound may be large. For the standard spherical harmonic basis on $S^2$ one has $K_\Lambda^2\asymp\Lambda^{\frac{1}{2}}$ as $\Lambda\to\infty$, or equivalently $K_{L(L+1)}^2\asymp L$ for maximum degree $L$, where the lower bound comes from the zonal harmonics and the upper bound from the addition formula. Weighted sampling and preconditioning can improve this dependence on the sphere and on certain surfaces of revolution, as shown by Burq, Dyatlov, Ward, and Zworski \cite{BDWZ}, whereas the theorem above uses uniform sampling and the chosen eigenbasis.
\end{remark}

\section{Concluding remarks}

All the estimates in this paper come from the same energy calculation. On $\mathbb Z_N^d$, the dimension and the number of controlled differences determine whether the inverse frequency weights are summable, logarithmically divergent, or polynomially divergent. In dimension two, one derivative is critical and leaves a square-root logarithm, for H\"older exponents below one the differences chosen separately on each frequency annulus remove the logarithmic loss, and with two derivatives the bound is independent of the grid size.

A common translation allows us to compare the discrete difference energy with the continuous energy without imposing a lower bound on the grid size, and averaging the zero coefficient as well retains the continuous mean. The translation that we select may depend on the function, so the problem of finding a favorable grid from limited observations is a separate question. On a compact surface, the coefficients $\frac{1}{\lambda_j}$ show that the critical square-root logarithm is necessary when the right-hand side is measured by the $L^2$ norm of the gradient, but additional regularity may improve the Fourier ratio, as the estimate with two derivatives on the square already shows.

Separate sensor displacements enter the reconstruction theorems as measurement noise. On the square we recover the sampled vector, and the piecewise constant extension controls the $L^2$ error on the torus, while on a manifold we reconstruct a function in a prescribed finite spectral space. In both cases the error bound consists of a contribution from spectral approximation and a contribution from placement accuracy. For uniform point sampling on a general manifold, the suprema of the chosen eigenfunctions remain in the sample count. Weighted sampling or a different eigenbasis may reduce this dependence, but any such change must also keep track of the Fourier ratio in the chosen basis.

\section*{Acknowledgments}

The authors gratefully acknowledge support from Williams College and the National Science Foundation under Grant No.~DMS-2241623 through the SMALL Undergraduate Research Program at Williams College.

\end{document}